\documentclass{amsart}
\usepackage{amssymb ,amsfonts}
\usepackage{xypic, color}
\usepackage{amsthm}

\newtheorem{thm}{Theorem}[section]
\newtheorem{cor}[thm]{Corollary}
\newtheorem{pro}[thm]{Proposition}
\newtheorem{lem}[thm]{Lemma}

\theoremstyle{definition}
\newtheorem{defn}[thm]{Definition}
\newtheorem{axiom}[thm]{Axiom}
\newtheorem{reformulation}[thm]{Reformulation}

\newtheorem{exmp}[thm]{Example}

\newtheorem{rem}[thm]{Remark}

\newcommand{\p}{\text{pro-}}

\newcommand{\ind}{\text{ind-}}
\DeclareMathOperator{\Hom}{Hom} \DeclareMathOperator{\Map}{Map}

\newcommand{\cc}{\mathcal{C}}

\newcommand{\rarr}{\rightarrow}
\newcommand{\map}{\rightarrow}

\newcommand{\ess}{\text{essentially levelwise }}
\DeclareMathOperator*{\colim}{colim}

\newcommand{\nn}{\mathbb{N}}
\newcommand{\zz}{\mathbb{Z}}

\newcommand{\inj}{\text{inj-}}
\newcommand{\proj}{\text{proj-}}

\newcommand{\Ar}{\text{Ar}}

\makeatletter
\let\c@equation\c@thm
\makeatother \numberwithin{equation}{section}

\newcommand{\dfn}{\textbf} 
\newcommand{\mdfn}[1]{\dfn{\mathversion{bold}#1}} 

\title{Model structures on pro-categories}
\author{Halvard Fausk}
\author{Daniel C. Isaksen}

\address{Department of  Mathematics\\ University of Oslo\\
1053 Blindern, 0316 Oslo, Norway}

\address{Department of Mathematics\\ Wayne State University\\
Detroit, MI 48202, USA}

\email{fausk@math.uio.no} \email{isaksen@math.wayne.edu}

\date{July 14, 2005}

\begin{document}

\begin{abstract}
We introduce a notion of a filtered model structure and use this
notion to produce various model structures on pro-categories.  This
framework generalizes the examples of \cite{pro}, \cite{strict}, and
\cite{ipi}. We give several examples, including a homotopy theory
for $G$-spaces, where $G$ is a pro-finite group.  The class of weak
equivalences is an approximation to the class of underlying weak
equivalences.
\end{abstract}

\maketitle

\section{Introduction}

The goal of this paper is to give a general framework for
constructing model structures on pro-categories. Given a proper
model structure on $\cc$, there is a strict model structure on $\p
\cc$ \cite{strict}. However, for most purposes, the class of strict
weak equivalences on pro-categories is too small.  For example,
generalized cohomology theories in the strict model structure on
pro-spaces do not typically have good computational properties. A
more useful class of maps are those that induce pro-isomorphisms on
pro-homotopy groups \cite{pro}. It turns out that this class of
equivalences is exactly the class of maps that are isomorphic to
strict $n$-equivalences for all integers $n$. In this paper we
axiomatize this situation so that it includes many other
interesting model structures on pro-categories. Also, we desire to
streamline the technical arguments that are usually required in
establishing a model structure on a pro-category. In the example
described above, the approach in this paper avoids many technical
issues from \cite{pro} involving basepoints.

The key idea is the notion of a {\em filtered} model structure. Our
main theorem (see Section \ref{sctn:ms}) states that every proper
filtered model structure on $\cc$ gives rise to a model structure on
$\p \cc$.

Because the list of axioms for a filtered model structure is
complicated (see Section \ref{sctn:filtered}), for now we will give
an example that gives a feeling for filtered model structures. Let
$\cc$ be a category, and let $A$ be a directed set such that $(W_a,
C_a, F_a)$ is a model structure on $\cc$ for every $a$ in $A$.
Moreover, assume that $W_a$ and $C_a$ are contained in $W_b$ and
$C_b$ respectively for $a \geq b$ (i.e., the classes $W_a$ and $C_a$
are ``decreasing'', so the classes $F_a$ are automatically
``increasing''). This is a particular example of a filtered model
structure, so there results an associated model structure on $\p
\cc$.  In $\p \cc$, a pro-map $f$ is a weak equivalence if for all
$a$ in $A$, $f$ is isomorphic to a pro-map that belongs to $W_a$
levelwise. The cofibrations are defined analogously.  The
fibrations, as usual, can be defined via a lifting property, but we
will give a more concrete description of them.

A filtered model structure is a generalization of the situation in the
previous paragraph.  We still have a directed set $A$ and
classes $W_a$, $C_a$, and $F_a$ of maps for each $a$ in $A$.
However, we do not assume that $(W_a, C_a, F_a)$ is necessarily a
model structure. For example, instead of requiring the
two-out-of-three property for each class $W_a$, we only require an
``up-to-refinement'' property.  Namely, for every $a$ in $A$, there
must exist a $b$ in $A$ such that if two out of the three maps $f$,
$g$, and $gf$ are in $W_b$, then the third is in $W_a$.  One
concrete example of this kind of phenomenon occurs when $A$ is the
set of natural numbers and $W_n$ is the class of $n$-equivalences of
spaces. This may seem like an unnatural generalization of the more
natural situation from the previous paragraph, but it is important
in our main examples.  The very nature of pro-categories suggests
that up-to-refinement definitions are a sensible approach.

We suspect that this axiomatization is more complicated than it has
to be, but we do not know any way to simplify it so that it still
includes all the examples of interest.

A t-model structure is a stable model structure
$\cc$ with a t-structure on its (triangulated) homotopy
category, together with a lift of the t-structure to $\cc$. This
notion is studied in detail in \cite{t-str}, where it is shown that
a particularly well-behaved filtered model structure on $\cc$
(and thus a model structure on $\p \cc$) can be associated to any
t-model structure.

This paper grew out of an attempt to find  useful model structures
on the category of pro-$G$-spaces and pro-$G$-spectra when $G$ is a
pro-finite group.  Section \ref{section-free} contains a detailed
description of a model structure on $G$-spaces as an illustration of
our general theory.
Analogous results for
pro-$G$-spectra are presented in detail in \cite{t-str} and
\cite{sfi}.

We now summarize our interest in pro-$G$-spaces when $G$ is
a pro-finite group.

Let $G$ be a finite group.  There is an obvious generalization to
pro-$G$-spaces of
the model structure for pro-non-equivariant spaces in the first paragraph.
Now
the weak equivalences are maps $f$ such that for every $n$, $f$ is
equivalent to a level map $g$ with the property that $g^H$ is a
levelwise $n$-equivalence for every subgroup $H$ of $G$. One can
make similar model structures for arbitrary topological groups.

When $G$ is a pro-finite group, the model structure on
pro-$G$-spaces described above is probably not the right
construction. For a pro-finite group, it is the continuous
cohomology $\colim_U H^*(G/U; M)$ rather than the group cohomology
$H^*(G^\delta; M)$ that is of interest.  Here $U$ ranges over the
finite-index normal subgroups of $G$, $M$ is a discrete continous
$G$-module, and $G^\delta$ is the group $G$ considered as a discrete
group. For finite groups, the group cohomology of $G$ is equal to
the Borel cohomology of a point.  In other words, the homotopy-orbit
space $*_{hG}$ of a point is $BG$.  In model-theoretic terms, what
is happening is that one takes a cofibrant replacement $EG$ for $*$
and then takes actual $G$-orbits to obtain $BG$.

We desire a model structure on pro-$G$-spaces such that the system
$\{ E(G/U) \}$ plays the role of the free contractible space $EG$.
In other words, $\{ E(G/U) \}$ should be a cofibrant replacement for
$*$.  Then $*_{hG}$ is equal to $\{ B(G/U) \}$, and the cohomology
of $*_{hG}$ is the continuous cohomology of $G$. Such a model
structure is obtained using our machinery by defining a weak
equivalence of pro-$G$-spaces to be a map such that for every
integer $n$, the map is isomorphic to a levelwise map $f$ with $f^U$
an $n$-equivalence for some finite-index subgroup $U$ of $G$.
The details of this particular example are given in
Section \ref{section-free}.

In fact, the $n$-equivalences in the previous paragraph are irrelevant
for the purposes of ordinary continuous cohomology.  We could just as well
consider a weak equivalence to be a map that is
isomorphic to a levelwise map $f$ with $f^U$
a weak equivalence for some finite-index subgroup $U$ of $G$.
The resulting model structure still behaves well with respect to
continuous cohomology, but it does not behave well with respect to
generalized continuous cohomology theories.  The non-equivariant analogue
of this phenomenon is explained in detail in \cite{ipi}.

In the context of equivariant model
categories that behave well with respect to continuous cohomology,
the paper \cite{goe} should also be mentioned.

\subsection{Organization}
We begin with a review of pro-categories, including a technical
discussion of essentially levelwise properties. Afterwards, we define
filtered model structures and prove our main  result in Theorem
\ref{thm:ms}, which establishes the existence of model structures on
pro-categories.  The next section considers Quillen functors in 
this context.

Then we proceed to examples. We first give a few  examples of our
general theory.   Then we focus on constructing  $G$-equivariant
homotopy theories when $G$ is a  pro-finite group.

\subsection{Background}
We assume that the reader is familiar with model categories,
especially in the context of equivariant homotopy theory.  The
original reference is \cite{Q}, but we will refer to more modern
treatments \cite{Hi} \cite{Ho}.  This paper is a generalization of
\cite{strict}, and we use specific pro-model category techniques
from it.

\section{Pro-categories}
We begin with a review of the necessary background on
pro-categories. This material can be found in \cite{SGA}, \cite{AM},
\cite{Duskin}, \cite{EH}, and \cite{strict}.

\subsection{Pro-Categories}

\begin{defn}
\label{defn:pro} For a category $\cc$, the category \mdfn{$\p \cc$}
has objects all cofiltering diagrams in $\cc$, and
$$\Hom_{\p \cc}(X,Y) = \lim_s \colim_t \Hom_{\cc}
     (X_t, Y_s).$$
Composition is defined in the natural way.
\end{defn}

A category $I$ is \dfn{cofiltering} if the following conditions
hold: it is non-empty and small; for every pair of objects $s$ and
$t$ in $I$, there exists an object $u$ together with maps $u \map s$
and $u \map t$; and for every pair of morphisms $f$ and $g$ with the
same source and target, there exists a morphism $h$ such that $fh$
equals $gh$. Recall that a category is \dfn{small} if it has only a
set of objects and a set of morphisms. A diagram is said to be
\dfn{cofiltering} if its indexing category is so. Beware that some
papers on pro-categories, such as \cite{AM} and \cite{Meyer},
consider cofiltering categories that are not small. All of our
pro-objects will be indexed by small categories.

Objects of $\p \cc$ are functors from cofiltering categories to
$\cc$. We use both set-theoretic and categorical language to discuss
indexing categories; hence ``$t \geq s$'' and ``$t \map s$'' mean
the same thing when the indexing category is actually a cofiltering
partially ordered set.

The word {\em pro-object} refers to an object of a pro-category. A
\dfn{constant} pro-object is one indexed by the category with one
object and one (identity) map. Let $\mathbf{c}: \cc \map \p \cc$ be the
functor taking an object $X$ to the constant pro-object with value
$X$. Note that this functor makes $\cc$ into a full subcategory of
$\p \cc$.

\subsection{Level Maps}

A {\bf level map} $X \map Y$ is a pro-map that is given by a natural
transformation (so $X$ and $Y$ must have the same indexing
category); this is a very special kind of pro-map. Up to
pro-isomorphism, every map is a level map \cite[App.~3.2]{AM}.

Let $M$ be a collection of maps in a category $\cc$. A level map $g$
in $\p \cc$ is a \mdfn{levelwise $M$-map} if  each $g_s$ belongs to
$M$. A pro-map is an \mdfn{essentially levelwise $M$-map} if it is
pro-isomorphic, in the category of arrows  in $\p \cc$, to a
levelwise $M$-map.

We will return to level maps in more detail in Section \ref{sctn:level}.

\subsection{Cofiniteness}

A partially ordered set $(I,\leq)$ is \dfn{directed} if for every
$s$ and $t$ in $I$, there exists $u$ such that $u \geq s$ and $u
\geq t$. A directed set $(I,\leq)$ is \dfn{cofinite} if for every
$t$, the set of elements $s$ of $I$ such that $s \leq t$ is finite.
A pro-object or level map is \dfn{cofinite directed} if it is
indexed by a cofinite directed set.

Every pro-object is isomorphic to a cofinite directed pro-object
\cite[Th.~2.1.6]{EH} (or \cite[Expos\'e~1, 8.1.6]{SGA}). Similarly,
up to isomorphism, every map is a cofinite directed level map.
Cofiniteness is critical for inductive arguments.

Let $f: X \map Y$ be a cofinite directed level map. For every index
$t$, the \mdfn{relative matching map $M_{t} f$} is the map
\[
X_{t} \map \lim_{s<t} X_{s} \times_{\lim_{s<t} Y_{s}} Y_{t}.
\]
The terminology is motivated by the fact that these maps appear in
Reedy model structures \cite[Defn.~16.3.2]{Hi}.

\begin{defn}
\label{defn:special}
Let $M$ be any class of maps in $\cc$. A map in $\p \cc$ is a
\mdfn{special $M$-map} if it is isomorphic to a cofinite directed
level map $f$ with the property that for each $t$, the relative matching map
$M_t f$ belongs to $M$.
\end{defn}

We will need the following lemma in several places later.  Its proof
is contained in the proof of \cite[Lem.~4.4]{strict}.

\begin{lem}
\label{lem:special} Let $f: X \map Y$ be a cofinite directed level
map.  For every $s$, the map $\lim_{t < s} X_t \map \lim_{t < s}
Y_t$ is a finite composition of base changes of the relative
matching maps $M_t f$ for $t < s$.
\end{lem}

\subsection{Simplicial structures on pro-categories}
\label{subsctn:simplicial}

Recall that a category $\cc$ is simplicial if for every object $X$
of $\cc$ and every simplicial set $K$, there exist objects $X
\otimes K$ (called a tensor) and $X^K$ (called a cotensor) of $\cc$
satisfying certain adjointness properties \cite[Sec.~9.1]{Hi}.
Moreover, for every $X$ and $Y$ in $\cc$, there is a simplicial set
$\Map(X,Y)$ that interacts appropriately with the tensor and
cotensor.

If $\cc$ is a simplicial category, then $\p \cc$ is again a
simplicial category.  Tensors and cotensors with finite simplicial
sets are defined levelwise, while tensors and cotensors with
arbitrary simplicial sets are defined via limits and colimits.  If
$X$ and $Y$ belong to $\p \cc$, then $\Map(X,Y)$ is defined to be
$\lim_s \colim_t \Map(X_t, Y_s)$.  See \cite[Sec.~16]{pro} for more
details.

\section{Essentially levelwise properties}
\label{sctn:level}

Later we will frequently encounter situations where a single pro-map
is an essentially levelwise $M_a$-map for all $a$, where $\{ M_a \}$ is a
collection of classes of maps
(see especially
Definitions \ref{defn:cofib}, \ref{defn:we}, and
\ref{defn:equiv-ms}).
This situation has some subtleties
that are worth exploring.

If $f$ is an essentially levelwise $M_a$-map for all $a$, it does not
follow that $f$ is an essentially levelwise $(\cap M_a)$-map.  The
problem is that different values of $a$ might require different
level maps, even though they are all pro-isomorphic to $f$.
However, we will show below in Corollary \ref{cor:ess-lev} 
that $f$ does have a slightly more
complicated but still concrete property.

We first need the following technical lemma for constructing 
isomorphisms in pro-categories.

\begin{lem}
\label{lem:pro-isomorphism}
Let $Y$ be a pro-object.
Suppose that for some of the maps $t \map s$ in the indexing diagram
for $Y$, there exists an object $Z_{ts}$ and a factorization
$Y_t \map Z_{ts} \map Y_s$ of the structure map $Y_t \map Y_s$.
Also suppose that for every $s$, there exists at least one $t \map s$
with this property.
The objects $Z_{ts}$ assemble into a pro-object $Z$ that is isomorphic
to $Y$.
\end{lem}

\begin{proof}
We may assume that $Y$ is indexed by a directed set $I$ 
because every pro-object is isomorphic to a pro-object indexed by
a directed set \cite[Thm.~2.1.6]{EH}.
Define a new directed set $K$ as follows.  The elements of 
$K$ consist of pairs $(t,s)$ of elements of $I$ such that
$t \geq s$ and a factorization $Y_t \map Z_{ts} \map Y_s$ exists.
If $(t,s)$ and $(t',s')$ are two elements of $K$, we say that
$(t',s') \geq (t,s)$ if $s' \geq t$.  It can easily
be checked that this makes $K$ into a directed set.

Note that the function $K \map I: (t,s) \mapsto s$
is cofinal in the sense of \cite[App.~1]{AM}.  
This means that we may reindex
$Y$ along this functor and assume that $Y$ is indexed by $K$;
thus we write $Y_{(t,s)} = Y_s$.

We define the pro-object $Z$ to be indexed by $K$ by setting
$Z_{(t,s)} = Z_{ts}$.  If $(t',s') \geq (t,s)$, then the structure map
$Z_{(t',s')} \map Z_{(t,s)}$ is the composition
\[
Z_{t's'} \map Y_{s'} \map Y_{t} \map Z_{ts},
\]
It can easily be checked that this gives a functor defined on $K$;
here is where we use that the composition $Y_t \map Z_{ts} \map Y_s$
equals $Y_t \map Y_s$.

Finally, we must show that $Z$ is isomorphic to $Y$.  We use the 
criterion from \cite[Lem.~2.3]{pro} for detecting pro-isomorphisms.
Given any $(t,s)$ in $K$, choose $u$ such that $(u,t)$ is in $K$.
Then there exists a diagram
\[
\xymatrix{
Z_{(u,t)} \ar[d]\ar[r] & Y_{(u,t)} = Y_t \ar[d]\ar[dl] \\
Z_{(t,s)} \ar[r] & Y_{(t,s)} = Y_s.                     }
\]
\end{proof}

In the rest of this section, we will frequently consider a collection
of classes $C_a$ indexed by a directed set such that $C_b$ is contained
in $C_a$ whenever $b \geq a$.  We say that the classes $C_a$ are a
\mdfn{decreasing collection} if they satisfy this property.

\begin{lem}
\label{lem:ess-lev}
Let $\{ C_a \}$ be a decreasing collection of classes of objects
indexed by a directed set $A$.  For each $a$ in $A$, let $C_a$ be a
class of objects such that $C_b$ is contained in $C_a$ when $b \geq a$.
An object in $\p \cc$ belongs to $C_a$ essentially levelwise for all $a$
if and only if it is isomorphic to a pro-object $X$ indexed by a
directed set $I$ such that for all $a$ in $A$, there exists an
element $s$ in $I$ (depending on $a$) with the property that $X_t$
belongs to $C_a$ for all $t \geq s$.
\end{lem}

The idea is that for every $a$, $X$
``eventually'' belongs to $C_a$ levelwise. However, the height at which
$X$ belongs to $C_a$ may depend on $a$.

\begin{proof}
Suppose that $X$ is a pro-object indexed by a
directed set $I$ such that for all $a$ in $A$, there exists an
element $s$ in $I$ (depending on $a$) with the property that $X_t$
belongs to $C_a$ for all $t \geq s$.
Given $a$,
choose $s$ such that $X_t$ belongs to $C_a$ for all $t \geq s$. Let
$I'$ be the subset of $I$ consisting of all $t \geq s$. By
restricting to $I'$, we obtain a new pro-object $X'$ such that $X'$
is isomorphic to $X$. By construction, $X'$ belongs to $C_a$
levelwise.
Since the above argument works for all $a$ in $A$, this
shows that $X$ belongs to $C_a$ essentially levelwise for every $a$.

Now we consider the other direction.
Let $X$ be a pro-object such that 
$X$ belongs to $C_a$ essentially levelwise for all $a$.
We may assume that $X$ is indexed by a directed set $I$ 
because every pro-object is isomorphic to a pro-object indexed by
a directed set \cite[Thm.~2.1.6]{EH}.
Reindex $X$ so that its indexing set is $I \times A$,
where $X_{s,a}$ equals $X_s$.  This only changes $X$ up to isomorphism.

For every  element $a$ of $A$,
choose an isomorphism $X \map Y^a$ where $Y^a$ is a pro-object that
belongs to $C_a$ levelwise.
The existence of this isomorphism implies that for any $s$ in $I$,
there exists $s' \geq s$ such that the structure map
$X_{s',a} \rarr X_{s,a}$ factors through $Y^a_{\phi(s,a)}$ for some 
$\phi(s,a)$ belonging to the indexing category of $Y^a$
(see \cite[Lem.~2.3]{pro} for a similar situation).
By Lemma \ref{lem:pro-isomorphism}, the objects $Y^a_{\phi(s,a)}$ 
assemble into a pro-object $Z$ indexed by $I \times A$ 
that is isomorphic to $X$.
By construction,
$Z$ has the property that $Z_{s,a'}$ belongs to $C_a$
for all $s$ in $I$ and all $a' \geq a$.
\end{proof}

\begin{cor}
\label{cor:ess-lev}
Let $\{ M_a \}$ be a decreasing collection of classes of maps indexed
by a directed set $A$.
A map in $\p \cc$ is an essentially levelwise $M_a$-map for all $a$
if and only if it has a level representation $f$ indexed by a
directed set $I$ such that for all $a$ in $A$, there exists an
element $s$ in $I$ (depending on $a$) with the property that $f_t$
belongs to $M_a$ for all $t \geq s$.
\end{cor}

\begin{proof}
This follows immediately from Lemma \ref{lem:ess-lev}, using that 
the category $\Ar(\p \cc)$ of arrows in $\p \cc$ is equivalent to the 
category $\p \Ar(\cc)$ \cite{Meyer}.
\end{proof}

We now discuss how functors preserve essentially levelwise properties.
This issue will be particularly relevant later in Section \ref{sctn:Quillen}
when we consider Quillen functors between pro-categories.

If $F: \cc \map \cc'$ is a functor, then $F$ induces another functor
$\p \cc \map \p \cc'$ in a natural way; we abuse notation and use
the same symbol $F$ for this functor. If $X = \{ X_s \}$ is a
filtered diagram, then $F X$ is defined to be the filtered diagram
$\{ F X_s \}$. 

\begin{lem} 
\label{lem:preserve} 
Let 
$\{ C_a \}$ and $\{ C'_{a'} \}$ be decreasing collections of classes
of objects in $\cc$ and
$\cc'$ respectively indexed on directed sets $A $ and $A'$ respectively.
Suppose  for every $ a' \in A' $ that
there is an $a \in A$ such that $F (C_a )$ is contained in $C_{a'}$.  
Consider the class of objects in $\p \cc$
that belong to $C_a$ essentially levelwise
for all $a$ in $A$.  
Then $F : \p \cc \rarr \p \cc'$ takes this class to the class of 
objects in $\p \cc'$ that belong to $C'_{a'}$ essentially levelwise
for all $a'$ in $A'$.
\end{lem} 

\begin{proof} 
For every $a'$, there is an $a$ such
that $F$ takes pro-objects that belong levelwise to $C_a$
to pro-objects that belong levelwise to $C'_{a'}$.
Finally, just recall that the functor
$F$ preserves isomorphisms.
\end{proof}

\begin{cor} 
\label{cor:preserve} 
Let 
$\{ M_a \}$ and $\{ M'_{a'} \}$ be decreasing collections of classes
of maps in $\cc$ and
$\cc'$ respectively indexed on directed sets $A $ and $A'$ respectively.
Suppose for every $ a' \in A' $ that
there is an $a \in A$ such that $F (M_a )$ is contained in $M'_{a'}$.  
Consider the class of maps in $\p \cc$
that are essentially levelwise $M_a$-maps for all $a$ in $A$.  
Then $F : \p \cc \rarr \p \cc'$ takes this class to the class of 
maps in $\p \cc'$ that are essentially levelwise $M'_{a'}$-maps 
for all $a'$ in $A'$.
\end{cor} 

\begin{proof}
This follows immediately from Lemma \ref{lem:preserve}, again using that
the categories $\Ar(\p \cc)$ and $\p \Ar(\cc)$ are equivalent.
\end{proof}

With  some assumptions, the sufficient condition given in Lemma
\ref{lem:preserve} is in fact necessary. Before we can prove this,
we need another technical lemma.

\begin{lem} 
\label{retract}
Let $\cc$ be a category and let $C$ be
a class of objects in  $\cc$ closed under retract. Let $X$
be a pro-object indexed by a directed set $I$ such that for all $t \geq s$,
the structure map $X_t \rarr X_s$ has a left inverse
(so  $X_t$ is a retract of $X_s$). If $X$ belongs to $C$
essentially levelwise,
then there is an element $s \in I$ such that $X_t $ is in $C$
for all  $t \geq  s$.
\end{lem}

\begin{proof}
Let $Y$ belong to $C$ levelwise such that $X$ is isomorphic to $Y$.
The existence of this isomorphism implies that for every $u$,
there exists $s \geq u$ and $v$ such that the structure map
$X_s \map X_u$ factors through $Y_v$.
Since $X_s$ is a retract of $X_u$, it follows that $X_s$ is also
a retract of $Y_v$.  But $Y_u$ belongs to $C$, so $X_s$ also belongs
to $C$.  

For any $t \geq s$, $X_t$ is a retract of $X_s$, so $X_t$
also belongs to $C$.
\end{proof}

\begin{lem}
Let $\cc$ and $\cc'$ be two categories such that $\cc$ is
cocomplete possessing an initial and terminal object $*$.
Let $\{C_a \}$ and $\{ C'_{a'} \}$ be decreasing
collections of classes of objects in $\cc$ and $\cc'$ respectively that
are indexed by directed sets $A$ and $A'$.
Suppose that each $C_a$ is closed under
small coproducts and that each $C'_{a'}$ is closed under retracts.
Let $F: \cc \rarr \cc'$ be a functor such that the associated functor $F :
\p \cc \rarr \p \cc'$ has the property that 
if $X$ belongs to $C_a$ essentially levelwise for all $a$ in $A$, then
$F(X)$ belongs to $C'_{a'}$ essentially levelwise for all $a'$ in $A'$.
Then for every $a'$ in $A'$, there exists  an element $a$ of $A$ such
that $F(C_a)$ is contained in $C'_{a'}$.
\end{lem}

\begin{proof}
Let $a' $ be an element in $ A' $.
Suppose for contradiction that for every $a$ in $A$,
$F (C_a )$ is not contained in $C'_{a'}$.
For each $b$ in $A$, we can choose an object $X_b$ in $C_b$ such that
 $F X_b$ is not in $C'_{a'}$. 
We construct a pro-object $Y$
indexed on $A$  by letting  $Y_a$ be
\[   \amalg_{c \geq a}   X_c.   \]
The structure maps are given by the canonical inclusions whenever
$b \geq a$.

Since the classes $C_a$ are decreasing and closed under coproducts,
each $Y_b$ belongs to $C_a$ for $b \geq a$.
Lemma \ref{lem:ess-lev} implies that $Y$ belongs to $C_a$ essentially
levelwise for all $a$.
Therefore, $FY$ belongs to $C'_{a'}$ essentially levelwise for all $a'$.

Since $\cc$ is pointed, each structure map $Y_b \map Y_a$ has a left inverse
given by projecting some of the factors to $*$.
Therefore each structure map $FY_b \map FY_a$ of $FY$ has a left inverse.
Lemma \ref{retract} applies, so for every $a'$, there exists
$b$ such that $FY_c$ belongs to $C'_{a'}$ whenever $c \geq b$.

Note that $X_a$ is a retract of $Y_a$ for all $a$, again using that $\cc$
is pointed.  Therefore $FX_a$ is a retract of $FY_a$ for all $a$.
We have shown in the previous paragraph that $FY_c$ belongs to $C'_{a'}$
for some value of $c$.  Since $C'_{a'}$ is closed under retract,
$FX_c$ also belongs to $C'_{a'}$.  This contradicts the way in which
we chose $X_c$.
\end{proof}

\begin{cor}
Let $\cc$ and $\cc'$ be two categories such that $\cc$ is
cocomplete possessing an initial and terminal object $*$.
Let $\{M_a \}$ and $\{ M'_{a'} \}$ be decreasing
collections of classes of maps in $\cc$ and $\cc'$ respectively that
are indexed by directed sets $A$ and $A'$.
Suppose that each $M_a$ is closed under
small coproducts and that each $M'_{a'}$ is closed under retracts.
Let $F: \cc \rarr \cc'$ be a functor such that the associated functor $F :
\p \cc \rarr \p \cc'$ has the property that 
if $f$ is an essentially levelwise $M_a$-map for all $a$ in $A$, then
$Ff$ is an essentially levelwise $M'_{a'}$-map for all $a'$ in $A'$.
Then for every $a'$ in $A'$, there exists  an element $a$ of $A$ such
that $F(M_a)$ is contained in $M'_{a'}$.
\end{cor}

\begin{proof}
This follows immediately from the previous lemma, again using that
the categories $\Ar(\p \cc)$ and $\p \Ar(\cc)$ are equivalent.
\end{proof}

\section{Filtered model structures}
\label{sctn:filtered}

We now describe the axiomatic setup that we use to produce model
structures on pro-categories.

\begin{defn}
\label{defn:filtered} A \mdfn{filtered model structure} consists of
a complete and cocomplete category $\cc$ equipped with:
\begin{enumerate}
\item
a directed set $A$,
\item
for each $a$ in $A$, a class $C_a$ of maps in $\cc$ such that $C_b$
is contained in $C_a$ whenever $b \geq a$,
\item
for each $a$ in $A$, a class $W_a$ of maps in $\cc$ such that $W_b$
is contained in $W_a$ whenever $b \geq a$,
\item
for each $a$ in $A$, a class $F_a$ of maps in $\cc$ such that $F_a$
is contained in $F_b$ whenever $b \geq a$.
\end{enumerate}
These classes must satisfy Axioms \ref{ax:2/3} through
\ref{ax:lift}, which are described below.
\end{defn}

\begin{axiom}
\label{ax:2/3} For all $a$ in $A$, there exists $b$ in $A$ such that
if $f$ and $g$ are two composable morphisms of $\cc$ with any two of
the three maps $f$, $g$, and $gf$ in $W_b$, then the third is in
$W_a$.
\end{axiom}

\begin{axiom}
\label{ax:retract} For all $a$ in $A$, the classes
$W_a$, $C_a$, and $F_a$ are closed under retract. 
The classes $C_a$ and $F_a$ are closed under composition.
The class $C_a$ is closed under arbitrary cobase change, and
the class $F_a$ is closed under arbitrary base change.
\end{axiom}

Recall that if $M$ is any class of maps, then \mdfn{$\inj M$} is
the class of maps that have the right lifting property with respect
to every map in $M$.  Similarly, \mdfn{$\proj M$} is the class of
maps that have the left lifting property with respect to every map
in $M$.

\begin{axiom}
\label{ax:acyclic} For every $a$ in $A$, the class $\inj C_a$ is
contained in $W_a \cap F_a$, and the class $\proj F_a$ is contained
in $W_a\cap C_a $.
\end{axiom}

\begin{axiom}
\label{ax:factor} For every map $f$ in $\cc$ and every $a$ in $A$,
$f$ can be factored as $p i$, where $i$ belongs to $C_a$ and $p$
belongs to $\inj C_a$, or as $q j$, where $j$ belongs to $\proj F_a$
and $q$ belongs to $F_a$.
\end{axiom}

\begin{axiom}
\label{ax:lift} If $f$ belongs to $W_a$ for some $a$ in $A$, then
there exists $b \geq a$ such that $f$ factors as $pi$, where $i$
belongs to $\proj F_b$ and $p$ belongs to $\inj C_b$.
\end{axiom}

\begin{defn}
\label{defn:F-inj-C} Let \mdfn{$F$} be the union of the classes
$F_a$ for all $a$. Let \mdfn{$\inj C$} be the union of the classes
$\inj C_a$ for all $a$.
\end{defn}

We think of $\{ F_a \}$ as an increasing filtration on $F$. Since
the classes $C_a$ are decreasing, the classes $\inj C_a$ are
increasing. Thus, $\{ \inj C_a \}$ is an increasing filtration on
$\inj C$.

\begin{defn}
A filtered model structure is \mdfn{proper} if it satisfies the
following two additional axioms.
\end{defn}

\begin{axiom}
\label{ax:cobase} For every $a$ in $A$, cobase changes of maps in
$W_a$ along maps in $C_a$ are in $W_a$.
\end{axiom}

\begin{axiom}
\label{ax:base} For every $a$ in $A$, base changes of maps in $W_a$
along maps in $F$ are in $W_a$.
\end{axiom}

Axioms \ref{ax:cobase} and \ref{ax:base} are a kind of properness.
They turn out to be necessary in various technical arguments
concerning pro-objects.  It may appear at first glance that Axiom
\ref{ax:base} is significantly stronger than Axiom \ref{ax:cobase}.
However, this is not the case.  Recall that the classes $C_b$ are
decreasing.

\begin{defn}
\label{defn:simplicial} A filtered model structure $\cc$ is
\mdfn{simplicial} if $\cc$ is a simplicial category satisifying the
following additional axiom.
\end{defn}

\begin{axiom}
\label{ax:simplicial} \mbox{}

\begin{enumerate}
\item
If $j: K \map L$ is a cofibration of finite simplicial sets and $i:
A \map B$ belongs to $C_a$, then the pushout-product map
\[
f: A \otimes L \amalg_{A \otimes K} B \otimes K \map B \otimes L
\]
belongs to $C_a$.
\item
If in addition $j$ is a weak equivalence or $i$ belongs to $\proj
F_a$, then $f$ belongs to $\proj F_a$.
\end{enumerate}
\end{axiom}

Axiom \ref{ax:simplicial} can be reformulated in the following two
equivalent ways.  The usual arguments with adjoints establish the
equivalence.

\begin{reformulation}
\mbox{}

\begin{enumerate}
\item
If $j: K \map L$ is a cofibration of finite simplicial sets and $p:
X \map Y$ belongs to $F_a$, then the map
\[
f: X^L \map Y^L \times_{Y^K} X^K
\]
belongs to $F_a$.
\item
If in addition $j$ is a weak equivalence or $p$ belongs to $\inj
C_a$, then $f$ belongs to $\inj C_a$.
\end{enumerate}
\end{reformulation}

\begin{reformulation}
\mbox{}

\begin{enumerate}
\item
If $i: A \map B$ belongs to $C_a$ and $p: X \map Y$ belongs to
$F_a$, then the map
\[
f: \Map(B,X) \map \Map(A,X) \times_{\Map(A,Y)} \Map(B,Y)
\]
is a fibration of simplicial sets.
\item
If in addition $i$ belongs to $\proj F_a$ or $p$ belongs to $\inj
C_a$, then $f$ is an acyclic fibration.
\end{enumerate}
\end{reformulation}

\begin{rem}
The axioms for a filtered model structure are almost but not quite
symmetric; see for example the inclusion relations for the classes
$C_a$ and $F_a$. The reason for this asymmetry is that the
construction of $\p \cc$ from $\cc$ is not symmetric.  If we were
interested in producing model structures on the ind-category $\ind
\cc$, then we would need to dualize the notion of a filtered model
category.
\end{rem}

At this point, we prove one simple lemma about filtered model structures
that we will need later.

\begin{lem}
\label{lem:inj-C} For any two elements $a$ and $b$ of $A$, the class
$\inj C_a$ is contained in the class $W_b$.
\end{lem}

\begin{proof}
Since $A$ is directed, we may choose $c$ such that $c \geq a$ and $c
\geq b$.  Since $C_c$ is contained in $C_a$, it follows formally that
$\inj C_a$ is contained in $\inj C_c$.  Now Axiom \ref{ax:acyclic}
implies that $\inj C_c$ is contained in $W_c$, which is contained in $W_b$.
\end{proof}

As an illustration of the definitions in this section, we provide several
general situations that produce filtered model structures.

\begin{pro} \label{pro-strict}
Suppose that $A$ consists of only a single element. A filtered model
structure (resp., proper filtered model structure, simplicial filtered
model structure)
indexed by $A$ is the same as an ordinary model structure
(resp., proper model structure, simplicial model structure)
on $\cc$.
\end{pro}

\begin{proof}
It is easy to verify that an ordinary model structure (resp.,
ordinary proper model structure, ordinary simplicial model
structure) gives a filtered model structure where $A$ has only one
element.

For the converse, suppose given a filtered model structure where $A$
has only one element.  We write $C$, $W$, and $F$ for the three
classes given in the definition. The two-out-of-three and retract
axioms are immediate from Axioms \ref{ax:2/3} and \ref{ax:retract}.
The factorization axiom follows from Axioms \ref{ax:acyclic} and
\ref{ax:factor}.

The lifting axiom requires more explanation.  Given a map $p$ in $W
\cap F$, use Axiom \ref{ax:lift} to factor it as $qj$, where $j$
belongs to $\proj F$ and $q$ belongs to $\inj C$.  By the retract
argument, $p$ is a retract of $q$.  Since $\inj C$ is formally
closed under retracts, it follows that $p$ also belongs to $\inj C$.
This shows that maps in $C$ lift with respect to maps in $W \cap F$.
The proof of the other half of the lifting axiom is identical. This
finishes the proof of the first claim.

For the second claim, if the filtered model structure is proper,
then properness for the ordinary model structure is given by Axioms
\ref{ax:cobase} and \ref{ax:base}.

Finally, for the third claim, if the filtered model structure is
simplicial, then Axiom \ref{ax:simplicial} implies that the ordinary
model structure is simplicial.  Note that for formal reasons (see,
for example, \cite[Prop.~16.1]{pro}), it is enough to check the
axioms for a simplicial model structure only for finite simplicial
sets.
\end{proof}

\begin{pro}
\label{pro:actual-mc} Let $A$ be a directed  set.  For each $a$ in
$A$, let $(C_a, W_a, F_a)$ be a proper model structure on $\cc$ such
that $C_a$ and $W_a$ are contained in $C_b$ and $W_b$ respectively
when $a \geq b$. Then $(A, C, W, F)$ is a proper filtered model
structure on $\cc$.
\end{pro}

\begin{proof}
Using that $\inj C_a$ equals $W_a \cap F_a$ and $\proj F_a$ equals
$W_a \cap C_a$, the verification of Axioms \ref{ax:2/3} through
\ref{ax:base} follow immediately from basic properties of model
structures.
\end{proof}

\begin{pro}
\label{pr:cof-gen}
Let $\cc$ be a proper cofibrantly generated model
category with a set $I$ of
generating cofibrations and a set $J$ of generating acyclic
cofibrations. Let $A$ be a directed set, and for each $a$ in $A$,
let $I_a$  be a subset of $I$ such that $I_a$ is contained in $I_b$
if $a \geq b$.
For each $a$ in $A$, define $C_a$ to be the class of all cofibrations,
$F_a$ to be
$\inj (I_a \cup J)$, and $W_a$ to be
the maps that are the composition of a map in $\proj F_a$
followed by an acyclic fibration.
If each $W_a$ is closed under retract and
Axiom \ref{ax:2/3} is satisfied, then $(A, C, W, F)$ is a filtered model
structure.
\end{pro}

Later in Section \ref{exmp-gen} we
describe several concrete examples
of this situation.
Beware that the filtered model structures arising from this proposition
are not necessarily proper.

\begin{proof}
We need to show that Axioms \ref{ax:2/3} through \ref{ax:lift} are
satisfied. We have assumed that Axiom \ref{ax:2/3} holds.

For Axiom \ref{ax:retract}, we have assumed that each $W_a$ is
closed under retract.  The class $C_a$ is closed under retract
because the class of cofibrations in $\cc$ is closed under
retract.  The class $F_a$ is closed under retract because it is
defined by a right lifting property.

The first part of Axiom \ref{ax:acyclic} is satisfied because $W_a$
and $F_a$ both contain the acyclic fibrations of $\cc$. For the
second part, every map in $\proj F_a$ is a cofibration because $F_a$
contains the acyclic fibrations of $\cc$.  Therefore, $\proj F_a$ is
contained in $C_a$.  It remains to explain why $\proj F_a$ is
contained in $W_a$, but this follows immediately from the definition
of $W_a$.

The first factorizations required by Axiom \ref{ax:factor} are just
the usual factorizations in $\cc$ of maps into cofibrations followed
by acyclic fibrations. The second factorizations can be produced by
applying the small object argument to the set $I_a \cup J$.

Axiom \ref{ax:lift} is satisfied by definition of $W_a$.
\end{proof}

\subsection{A possible weakening of the axioms}
For expository reasons, we have chosen a form of 
Axiom \ref{ax:retract} that is unnecessarily strong. To develop most
of our theory,
it suffices to assume only that the classes $C_a$ and $W_a$
are closed under retract.
In any case, we can always
replace $C_a$ and $F_a$ by the smallest classes of maps
containing $C_a$ and $F_a$
respectively and satisfying the conditions in \ref{ax:retract}. These
new classes (with $W_a$ unchanged)
will still satisfy all of the rest of the axioms for a filtered model
structure.  However,
Axioms \ref{ax:cobase}, \ref{ax:base}, and \ref{ax:simplicial} might 
no longer be satisfied.

Axioms \ref{ax:cobase}, \ref{ax:base}, and \ref{ax:simplicial} are
not up-to-refinement properties, unlike Axiom \ref{ax:2/3}.  In
fact, there is an obvious way to weaken these three axioms to make
them up-to-refinement.  We have chosen to not formalize this in our
definitions because we know of no examples in which the added
generality is necessary.

The more general form of Axiom \ref{ax:cobase} states that
for every $a$ in $A$,
there exists $b$ and $c$ in $A$ such that cobase changes of maps in
$W_b$ along maps in $C_c$ are in $W_a$.

The more general form of Axiom \ref{ax:base} states that for every $a$ in $A$,
there exists $b$ in $A$ such that base changes of maps in $W_b$ along maps
in $F$ are in $W_a$.

The more general form of Axiom \ref{ax:simplicial} states that
for every $a$ in $A$, there exists $b$ in $A$ such that
if $j: K \map L$ is a cofibration of finite simplicial sets and $i:
A \map B$ belongs to $C_b$, then the pushout-product map
\[
f: A \otimes L \amalg_{A \otimes K} B \otimes K \map B \otimes L
\]
belongs to $C_a$;
if in addition $j$ is a weak equivalence or $i$ belongs to $\proj
F_b$, then $f$ belongs to $\proj F_a$.

\section{Model structures on $\p \cc$}
\label{sctn:ms}

Throughout this section, let $\cc$ be a category equipped with a
proper filtered model structure $(A, C, W, F)$ as in the
previous section.

\begin{defn}
\label{defn:cofib} A map in $\p \cc$ is a \mdfn{cofibration} if it
is an essentially levelwise $C_a$-map for every $a$ in $A$.
\end{defn}

\begin{defn}
\label{defn:we} A map in $\p \cc$ is a \mdfn{weak equivalence} if it
is an essentially levelwise $W_a$-map for all $a$ in $A$.
\end{defn}

If $i$ is a cofibration, there is no guarantee that $i$ has one
level replacement that is a level $C_a$-map for every $a$.
Typically, the choice of level replacement depends on $a$.  The same
warning applies to weak equivalences.

Corollary \ref{cor:ess-lev} can be used
to give a slightly more concrete description of the cofibrations and
weak equivalences.  For understanding specific examples, this more
concrete description is often helpful.  However, for proving general
results, we prefer to work with the more abstract definition.

\begin{defn}
\label{defn:fib} A map in $\p \cc$ is a \mdfn{fibration} if it is a
retract of a special $F$-map.
\end{defn}

Recall that an \mdfn{acyclic cofibration} is a map that is both a
cofibration and a weak equivalence.  Similarly, an \mdfn{acyclic
fibration} is a map that is both a fibration and a weak equivalence.

We will eventually prove that these definitions yield a model
structure on $\p \cc$. First we need a series of lemmas that will
lead to the proof. Our approach follows \cite{strict}.

\begin{lem} \label{lem:2/3}
Suppose that $f$ and $g$ are two composable morphisms in $\p \cc$.
If any two of $f$, $g$, and $gf$ are weak equivalences, then so is
the third.
\end{lem}

\begin{proof}
Suppose that two of $f$, $g$, and $gf$ are weak equivalences.  We
need to show that the third is an essentially levelwise $W_a$-map
for every $a$. Choose $b$ such that if any two of $\phi$, $\psi$,
and $\psi\phi$ are in $W_b$, then the third is in $W_a$; this is
possible by Axiom \ref{ax:2/3}.

If we assume that any two of $f$, $g$, and $gf$ are essentially
levelwise $W_b$-maps, then the proofs of \cite[Lem.~3.5]{strict} and
\cite[Lem.~3.6]{strict} can be applied to conclude that the third is
an essentially levelwise $W_a$-map. Note that
\cite[Lem.~3.2]{strict} works for $C_b$ and $\inj C_b$ (or for
$\proj F_b$ and $F_b$) because of Axiom \ref{ax:factor} (see
\cite[Rem.~3.3]{strict}). To make these proofs work, we need Axioms
\ref{ax:cobase} and \ref{ax:base}.

To illustrate the point, we describe in detail how to adapt the
proof of \cite[Lem.~3.5]{strict} to our situation.  Suppose that $f$
and $g$ are essentially levelwise $W_b$-maps.  We wish to show that
$gf$ is an essentially levelwise $W_a$-map.

We may assume that $f$ and $g$ are levelwise $W_b$-maps, but their
index categories are not necessarily the same.  However, we can
obtain a levelwise diagram
\[
\xymatrix{ X \ar[r]^f & Y & Z \ar[l]_h^{\cong} \ar[r]^g & W}
\]
in which $f$ and $g$ are levelwise $W_b$-maps while $h$ is a
pro-isomorphism (but not a levelwise isomorphism). We must construct
a levelwise $W_a$-map isomorphic to the composition $gh^{-1}f$.

By \cite[Lem.~3.2]{strict}, after reindexing we can factor $h: Z
\map Y$ into a levelwise $C_b$-map $Z \map A$ followed by a
levelwise $F_b$-map $A \map Y$ such that both maps are pro-isomorphisms.
Here we are using Axiom
\ref{ax:factor} to provide the necessary factorizations in $\cc$
(and also Axiom \ref{ax:acyclic} to identify that a map in $\inj
C_a$ is necessarily in $F_a$). We now have a diagram
\[
\xymatrix{ X \ar[r] & Y & A \ar[l]_{\cong} & Z \ar[l]_{\cong} \ar[r]
& W }
\]
in which the first and fourth maps are levelwise $W_b$-maps, and the
second and third are pro-isomorphisms.

Let $B$ be the pullback $X \times_Y A$, and let $C$ be the pushout
$A \amalg_Z W$, which we may construct levelwise. The map $B \map A$
is levelwise a base change of a map in $W_b$ along a map in $F_b$.
By Axiom \ref{ax:base}, $B \map A$ is a levelwise $W_b$-map.
Similarly, the map $A \map C$ is levelwise a cobase change of a map
in $W_b$ along a map in $C_b$.  By Axiom \ref{ax:cobase}, $A \map C$
is a levelwise $W_b$-map.

The maps $B \map X$ and $W \map C$ are pro-isomorphisms since base
and cobase changes preserve isomorphisms. Hence the composition $B
\map C$ is isomorphic to $gh^{-1}f$ as desired. Moreover, $B \map C$
is levelwise a composition of two maps in $W_b$, which means that it
is a levelwise $W_a$-map because of the way in which $b$ was chosen.
\end{proof}

Recall the following lemma from \cite[Thm.~5.5]{lim}.

\begin{lem}
\label{lem:retract} Let $M$ be any class of maps in $\cc$.  Then the
class of essentially levelwise $M$-maps in $\p \cc$ is closed under
retracts.
\end{lem}

\begin{cor}
\label{cor:retract} The class of cofibrations and the class of weak
equivalences in $\p \cc$ are closed under retract.
\end{cor}

\begin{proof}
The class of cofibrations is the intersection of a set of classes,
each of which is closed under retract by Lemma \ref{lem:retract}.
The same argument applies to the weak equivalences.
\end{proof}

\begin{lem}
\label{lem:factor1} Every map $f:X \map Y$ in $\p \cc$ factors as a
cofibration $i:X \map Z$ followed by a special $\inj C$-map $p:Z
\map Y$.
\end{lem}

\begin{proof}
We may suppose that $f$ is a level map indexed by a cofinite
directed set $I$.  Moreover, by adding isomorphisms to cofiltered
diagrams, we may assume that $I$ has cardinality larger than $A$.
Choose an arbitrary function $\phi: I \map A$ such that $\phi(s)
\geq \phi(t)$ if $s \geq t$ and such that for all $a$ in $A$, there
exists $s$ in $I$ such that $\phi(s) \geq a$; this function can be
constructed inductively because $I$ is cofinite and because the
cardinality of $I$ is larger than the cardinality of $A$.

Suppose for induction that the maps $i_t: X_t \map Z_t$ and $p_t:Z_t
\map Y_t$ have already been defined for $t < s$. Consider the map
\[
X_s \map Y_s \times_{\lim_{t<s} Y_t} \lim_{t<s} Z_t.
\]
Use Axiom \ref{ax:factor} to factor it into a map $i_s: X_s \map
Z_s$ belonging to $C_{\phi(s)}$ followed by a map
\[
q_s: Z_s \map Y_s \times_{\lim_{t<s} Y_t} \lim_{t<s} Z_t
\]
belonging to $\inj C_{\phi(s)}$. Let $p_s$ be the map $Z_s \map Y_s$
induced by $q_s$. This extends the factorization to level $s$.

It follows immediately from its construction that $p$ is a special
$\inj C$-map.  To show that $i$ is a cofibration, just apply Lemma
\ref{lem:ess-lev}.
\end{proof}

\begin{lem}
\label{lem:factor2} Every map $f:X \map Y$ in $\p \cc$ factors into
a map $i:X \map Z$ followed by a special $F$-map $p:Z \map Y$, where
$i$ is an essentially levelwise $\proj F_a$-map for every $a$.
\end{lem}

\begin{proof}
The proof is identical to the proof of Lemma \ref{lem:factor1},
except that we factor the map
\[
X_s \map Y_s \times_{\lim_{t<s} Y_t} \lim_{t<s} Z_t
\]
into a map belonging to $\proj F_{\phi(s)}$ followed by a map
belonging to $F_{\phi(s)}$.
\end{proof}

\begin{lem}
\label{lem:lift1} A map in $\p \cc$ is a cofibration if and only if
it has the left lifting property with respect to all retracts of
special $\inj C$-maps. Also, a map in $\p \cc$ is a retract of a
special $\inj C$-map if and only if it has the right lifting
property with respect to all cofibrations.
\end{lem}

\begin{proof}
First we will show that cofibrations have the left lifting property
with respect to retracts of special $\inj C$-maps.
Let $i:A \map B$ be a cofibration, and let $p$ be a retract of a
special $\inj C$-map. Since retracts preserve lifting properties, it
suffices to assume that $p$ is a special $\inj C$-map.  Moreover, as
shown in \cite[Prop.~5.2]{strict}, a special $\inj C$-map is a
composition along a transfinite tower, each of whose maps is a base
change of a map of the form $cX \map cY$, where $X \map Y$ belongs
to $\inj C$. Since base changes and transfinite compositions
preserve lifting properties, it suffices to assume that $p$ is the
map $cX \map cY$, where $X \map Y$ belongs to $\inj C$.

The map $X \map Y$ belongs to $\inj C_a$ for some $a$.  Since $i$ is
a cofibration, it is an essentially levelwise $C_a$-map.  Therefore,
we may assume that $i$ is a levelwise $C_a$-map.

Suppose given a square
\[
\xymatrix{
A \ar[r] \ar[d]_{i} & cX \ar[d]^p \\
B \ar[r] & cY                      }
\]
in $\p \cc$.  This square is represented by a square
\[
\xymatrix{
A_s \ar[r] \ar[d]_{i_s} & X \ar[d] \\
B_s \ar[r] & Y                      }
\]
in $\cc$.  This square has a lift because $X \map Y$ has the right
lifting property with respect to all maps in $C_a$.  Finally, this
lift represents a map $B \map cX$ that is our desired lift.

Now suppose that a map $i: A \map B$ has the left lifting property
with respect to all special $\inj C$-maps. Use Lemma
\ref{lem:factor1} to factor $i$ as a cofibration $i':A \map B'$
followed by a special $\inj C$-map $p:B' \map B$. Since $i$ has the
left lifting property with respect to $p$ by assumption, the retract
argument implies that $i$ is a retract of $i'$. But retracts
preserve cofibrations by Corollary \ref{cor:retract}, so $i$ is
again a cofibration.

 Finally, suppose that $p: X \map Y$
has the right lifting property with respect to all cofibrations. Use
Lemma \ref{lem:factor1} to factor $p$ as a cofibration $i:X \map X'$
followed by a special $\inj C$-map $p':X' \map Y$. Similarly to the
previous paragraph, $p$ is a retract of $p'$, so $p$ is a retract of
a special $\inj C$-map.
\end{proof}

\begin{lem}
\label{lem:lift2} A map in $\p \cc$
is an essentially levelwise $\proj F_a$-map for all $a$ if and only if it
has the left lifting property with respect to all fibrations.
Also, a map in $\p \cc$ is a
fibration if and only if it has the right lifting property with
respect to all essentially levelwise $\proj F_a$-maps for all $a$.
\end{lem}

\begin{proof}
The proof is the same as the proof of Lemma \ref{lem:lift1}, except
that the role of cofibrations is replaced by the maps that are
essentially levelwise $\proj F_a$-maps for every $a$ and special
$\inj C$-maps are replaced by special $F$-maps. Lemma
\ref{lem:factor2} is relevant instead of Lemma \ref{lem:factor1}.
Note also that retracts preserve the maps that are essentially
levelwise $\proj F_a$-maps for every $a$; the proof is like the
proof of Corollary \ref{cor:retract}.
\end{proof}

\begin{pro}
\label{pro:a-cofib} A map in $\p \cc$ is an acyclic cofibration if
and only if it is an essentially levelwise $\proj F_a$-map for every
$a$.
\end{pro}

\begin{proof}
One implication follows from the definitions and the fact that
$\proj F_a$ is contained in $W_a \cap C_a$ by Axiom
\ref{ax:acyclic}.

For the other implication, let $i:A \map B$ be a weak equivalence
and cofibration. Fix an arbitrary $a$; we will show that $i$ is an
essentially levelwise $\proj F_a$-map.

Given $a$, begin by choosing $b$ as in Axiom \ref{ax:2/3}.
We may assume that $i$ is a level map such that each $i_s$ belongs
to $W_b$. As in the proof of Lemma \ref{lem:factor1},
we may use Axiom \ref{ax:lift}
to factor $i$ into a map $i': A \map B'$ followed
by a map $p: B' \map B$ such that for each $s$, $i'_s$ belongs to
$\proj F_c$ and the relative matching map
\[
M_s p: B'_s \map B_s \times_{\lim_{t<s} B_t} \lim_{t<s} B'_t
\]
belongs to $\inj C_c$ for some $c \geq a$.

Actually, in order to apply Axiom \ref{ax:lift}, we need to prove
inductively that the map
$A_s \map B'_s \times_{\lim_{t<s} B'_t} \lim_{t<s} B_t$ belongs to
$W_a$.  To do this, consider the diagram
\[
\xymatrix{
A_s \ar[r]\ar[dr] & B_s \times_{\lim_{t<s} B_t} \lim_{t<s} B'_t \ar[r]\ar[d] &
    \lim_{t<s} B'_t \ar[d] \\
& B_s \ar[r] & \lim_{t<s} B_t.  }
\]
Using Axiom \ref{ax:2/3}, we just need to show that the diagonal
map and the left vertical map belong to $W_b$.  The diagonal map
belongs to $W_b$ by the assumption on $i$.
On the other hand, the right vertical
map belongs to $\inj C$ by Lemma \ref{lem:special} and the induction
assumption.  Therefore, the left vertical map is also in $\inj C$
since $\inj C$ is closed under base changes for formal reasons.
Now Lemma \ref{lem:inj-C} implies that it belongs to $W_b$.

At this point, we have factored $i$ as $pi'$, where
$p$ is a special $\inj C$-map.
Since $i$ is a cofibration, it has the left lifting
property with respect to $p$ by Lemma \ref{lem:lift1}. The retract
argument now implies that $i$ is a retract of $i'$. Because
essentially levelwise $\proj F_a$-maps are closed under retract by
Lemma \ref{lem:retract}, it suffices to show that $i'$ is a
levelwise $\proj F_a$-map. Recall that each $i'_s$ belongs to $\proj
F_c$ for some $c \geq a$. Since $F_a$ is contained in $F_c$, it
follows that $\proj F_c$ is contained in $\proj F_a$; therefore,
each $i'_s$ belongs to $\proj F_a$.
\end{proof}

The following proposition, although not actually necessary to
establish the existence of the model structure on $\p \cc$, is a
useful detection principle for acyclic cofibrations.  It says that
the fibrations are ``generated'' by a certain very simple class of
fibrations.

\begin{pro}
\label{pro:detect-acyclic-cofibrations} A map $i: A \map B$ is an
acyclic cofibration if and only if it has the left lifting property
with respect to all constant pro-maps $cX \map cY$ in which $X \map
Y$ belongs to $F$.
\end{pro}

\begin{proof}
Lemma \ref{lem:lift2} and Proposition \ref{pro:a-cofib} imply that
$i$ is an acyclic cofibration if and only if it has the left lifting
property with respect to all special $F$-maps. By
\cite[Prop.~5.2]{strict}, every special $F$-map is a transfinite
composition of a tower of maps, each of whose maps is a base change
of a map of the form $cX \map cY$ with $X \map Y$ in $F$.  By a
formal argument with lifting properties, a map has the left lifting
property with respect to all special $F$-maps if and only if it has
the left lifting property with respect to maps of the form $cX \map
cY$ with $X \map Y$ in $F$.
\end{proof}

\begin{pro} \label{pro:a-fib}
A map $p$ is an acyclic fibration if and only if it is a retract of
a special $\inj C$-map.
\end{pro}

\begin{proof}
First suppose that $p$ is a retract of a special $\inj C$-map. The
class of acyclic fibrations is closed under retract by Corollary
\ref{cor:retract} and by the definition of fibrations, so it
suffices to assume that $p$ is a special $\inj C$-map. Since each
class $\inj C_a$ is contained in $F_a$ by Axiom \ref{ax:acyclic}, it
follows that the union $\inj C$ is contained in the union $F$.
Therefore, every special $\inj C$-map is a special $F$-map.  This
shows that $p$ is a fibration.

It remains to show that $p$ is a weak equivalence. Given a fixed
$s$, we will show that $p_s: X_s \map Y_s$ belongs to $W_a$ for all
$a$.  Then $p$ is a level $W_a$-map for every $a$ and thus a weak
equivalence. In order to do this, Lemma \ref{lem:inj-C} tells us
that we only have to show that $p_s$ belongs to $\inj C$.

We may choose an element $a$ of $A$ such that $M_t p$ belongs to
$\inj C_a$ for every $t \leq s$.  This follows from cofiniteness and
the fact that the classes $\inj C_b$ are increasing.

The map $p_s: X_s \map Y_s$ factors as
\[
\xymatrix@1{ X_s \ar[r]^-{M_s p} & Y_s \times_{\lim_{t<s} Y_t}
\lim_{t<s} X_t \ar[r]^-{q_s} & Y_s. }
\]
Our goal is to show that $p_s$ belongs to $\inj C_a$. Since
compositions and base changes preserve $\inj C_a$ for formal
reasons, it suffices to show that $\lim_{t<s} p_t: \lim_{t<s} X_t
\map \lim_{t<s} Y_t$ belongs to $\inj C_a$. This last map is a
finite composition of maps that are base changes of the maps $M_t p$
for $t \leq s$ (see Lemma \ref{lem:special}).
Finally, recall that each $M_t p$ belongs
to $\inj C_a$. This finishes one implication.

Now suppose that $p:X \map Y$ is an acyclic fibration.  Use Lemma
\ref{lem:factor1} to factor $p$ into a cofibration $i: X \map X'$
followed by a special $\inj C$-map $p': X' \map Y$.  By the
two-out-of-three axiom (see Lemma \ref{lem:2/3}), we know that $i$
is in fact an acyclic cofibration.  Then Proposition
\ref{pro:a-cofib} says that $i$ is an essentially levelwise $\proj
F_a$-map for every $a$, so Lemma \ref{lem:lift2} implies that $p$
has the right lifting property with respect to $i$.  The retract
argument then gives that $p$ is a retract of $p'$, as desired.
\end{proof}

The following lemma will be needed to show that the model structure
on $\p \cc$ is right proper.

\begin{lem} \label{lem:fib}
Any special $F$-map is a levelwise $F$-map.
\end{lem}

\begin{proof}
Suppose given a cofinite directed level map $p: X \map Y$ for which
each relative matching map $M_s p$ belongs to $F$. The map $p_s$ is the
composition of $M_s p$ followed by the projection $Y_s
\times_{\lim_{t <s} Y_t} \lim_{t<s} X_t \map Y_s$. This projection
is a base change of the map $\lim_{t<s} X_t \map \lim_{t<s} Y_t$,
which is a finite composition of base changes of the maps $M_t p$
for $t < s$ by Lemma \ref{lem:special}. So Axiom \ref{ax:retract}
implies that $p_s$ belongs to $F$.
\end{proof}

\begin{thm}
\label{thm:ms} Let $ (A, C ,W, F)$ be a proper filtered model
structure on $\cc$. Then  $\p \cc$ has a proper model structure
given by Definitions \ref{defn:cofib}, \ref{defn:we}, and
\ref{defn:fib}.
\end{thm}

\begin{proof}
The category $\p \cc$ is complete and cocomplete because $\cc$ is
complete and cocomplete \cite[Prop.~11.1]{pro}. The two-out-of-three
axiom for weak equivalences is not automatic;  we proved this in
Lemma \ref{lem:2/3}. Corollary \ref{cor:retract} shows that
cofibrations and weak equivalences are closed under retract.
Fibrations are closed under retract by definition.

Factorizations into cofibrations followed by acyclic fibrations are
given in Lemma \ref{lem:factor1}.  Here we have to use Proposition
\ref{pro:a-fib} to identify the second map as an acyclic fibration.
Factorizations into acyclic cofibrations followed by fibrations are
given in Lemma \ref{lem:factor2}.  Now we have to use Proposition
\ref{pro:a-cofib} to identify the first map as an acyclic
cofibration.

Cofibrations lift with respect to acyclic fibrations by Lemma
\ref{lem:lift1}; we use Proposition \ref{pro:a-fib} to identify the
acyclic fibrations. Acyclic cofibrations lift with respect to
fibrations by Lemma \ref{lem:lift2}; now we use Proposition
\ref{pro:a-cofib} to identify the acyclic cofibrations.

We now show that the model structure is proper.
 For right
properness, consider a pullback square
\[
\xymatrix{
W \ar[r]^q \ar[d]_g & X \ar[d]^f \\
Y \ar[r]_p & Z               }
\]
in which $f$ is a weak equivalence and $p$ is a fibration. We want
to show that $g$ is also a weak equivalence; that is, we want to
show that $g$ is an essentially levelwise $W_a$-map for all $a$.

Lemma \ref{lem:fib} says  that $p$ is an essentially levelwise
$F$-map. Therefore, the proof of \cite[Thm.~4.13]{strict} can be
applied to show that base changes of essentially levelwise
$W_a$-equivalences along fibrations are essentially levelwise
$W_a$-equivalences. We need that $F$ is closed under arbitrary base
changes and Axiom \ref{ax:base} for the proof to work. Since $f$ is
an essentially levelwise $W_a$-equivalence, we conclude that $g$ is
an essentially levelwise $W_a$-equivalence.

The proof of left properness is dual but easier because we know from
the definition that cofibrations are essentially levelwise
$C_a$-maps. Now we need to use that $C_a$ is closed under arbitrary
cobase changes and Axiom \ref{ax:cobase}.
\end{proof}

\begin{thm}
\label{thm:simplicial} Let $ ( A, C ,W, F)$ be a simplicial proper
filtered model structure on $\cc$. Then the model structure of
Theorem \ref{thm:ms} is also simplicial.
\end{thm}

\begin{proof}
As explained in Section \ref{subsctn:simplicial}, $\p \cc$ is a
simplicial category. Let $j: K \map L$ be a cofibration of
simplicial sets, let $i: X \map Y$ be a cofibration in $\p \cc$, and
let $f$ be the map
\[
f: X \otimes L \amalg_{X \otimes K} Y \otimes K \map Y \otimes L.
\]
As explained in \cite[Prop.~16.1]{pro}, it suffices to assume that
$K$ and $L$ are finite simplicial sets.

Let $a$ be any element of $A$. We may assume that $i$ is a levelwise
$C_a$-map.  Because $K$ and $L$ are finite, the map $f$ may be
constructed levelwise; that is, $f_s$ equals
\[
X_s \otimes L \amalg_{X_s \otimes K} Y_s \otimes K \map Y_s \otimes
L.
\]
Because of part (1) of Axiom \ref{ax:simplicial}, $f$ is a levelwise
$C_a$-map. This shows that $f$ is a cofibration in $\p \cc$.

Next, assume that $j$ is an acyclic cofibration.  As in the previous
paragraph but using part (2) of Axiom \ref{ax:simplicial}, $f$ is an
essentially levelwise $\proj F_a$-map for every $a$.  By Proposition
\ref{pro:a-cofib}, it follows that $f$ is an acyclic cofibration.

Finally, assume that $i$ is an acyclic cofibration. By Proposition
\ref{pro:a-cofib}, $i$ is an
essentially levelwise $\proj F_a$-map for every $a$. As before but
using the other part of part (2) of Axiom \ref{ax:simplicial}, $f$
is an essentially levelwise $\proj F_a$-map for every $a$, so it is
an acyclic cofibration by Proposition \ref{pro:a-cofib}.
\end{proof}

The following result shows that the model structures produced
by Theorem \ref{thm:ms} are a generalization of
the strict model structures of \cite{EH} and \cite{strict}.

\begin{pro} \label{exmp:strict}
Let $\cc$ be a proper model category considered as a proper filtered
model category indexed on a set $A$ with only one element.
The associated model structure on $\p \cc$ is the strict
model structure on $\p \cc$. 
\end{pro}

\begin{proof}
This follows from the definitions, Theorem \ref{thm:ms}, and 
Proposition \ref{pro-strict}.
\end{proof}

Recall from \cite{strict} that if $\cc$ is simplicial, then the
strict model structure on $\p \cc$ is also simplicial.  This result
now is an immediate corollary of Theorem \ref{thm:simplicial}.

We include one more minor technical lemma in this section that will
be needed later.

\begin{lem}
\label{lem:level-fibrant}
Suppose that $Y$ is a fibrant object in $\p \cc$.  Then $Y$ is isomorphic
to an object $Y'$ such that each map $Y'_s \map *$ belongs to $F$.
\end{lem}

\begin{proof}
Let $D$ be the class of objects $X$ in $\cc$ such that $X \map *$ 
belongs to $F$.  We want to show that $Y$ belongs to $D$ essentially
levelwise.

We may assume that $Y$ is cofinite directed.  
If we use Lemma \ref{lem:factor2}
to factor $Y \map *$ into an acyclic cofibration
followed by a fibration and then apply the retract argument, we see that
$Y \map *$ is a retract of a special $F$-map $Z \map *$.
Lemma \ref{lem:fib} implies that $Z$ belongs to $D$ levelwise,
and \cite[Thm.~5.5]{lim} implies that $Y$ belongs to $D$
essentially levelwise.
\end{proof}

\section{Quillen functors}
\label{sctn:Quillen}

In this section, we consider two proper filtered model structures
$(A, C, W, F)$ and $(A', C', W', F')$ on the categories $\cc$ and
$\cc'$ respectively.   We will compare the associated model
structures on $\p \cc$ and $\p \cc'$.

Recall that if
if $L: \cc \map \cc'$ is a functor, then $L$ induces another functor
$\p \cc \map \p \cc'$ by applying $L$ levelwise.
Moreover, if $R: \cc' \map \cc$ is the right adjoint
of $L$, then $R: \p \cc' \map \p \cc$ is the right adjoint of $L: \p
\cc \map \p \cc'$.  

We will give precise
conditions telling us when the induced functors $L: \p \cc \map \p
\cc'$ and $R: \p \cc' \map \p \cc$ are a Quillen adjoint pair.
Then we will give some conditions that
imply that $L$ and $R$ induce a Quillen equivalence between
$\p \cc$ and $\p \cc'$.

\begin{thm}
\label{thm:Q-adj}
Let $(A, C, W, F)$ and $(A', C', W', F')$ be
proper filtered model structures on the categories $\cc$ and $\cc'$
respectively. Let $L : \cc \rarr \cc'$ be a left adjoint of $R :
\cc' \rarr \cc$.
The induced functors $L: \p \cc \rarr \p \cc'$ and 
$R : \p \cc' \rarr \p \cc$ are a Quillen
adjoint pair when $\p \cc$ and $\p \cc'$ are equipped with the model
structures of Theorem \ref{thm:ms} if and only if
$R(F')$ is contained in $F$ and
$R(\inj C')$ is contained in $\inj C$. 
\end{thm}

The properness assumption on the filtered model structures is only
to guarantee that the model structures on $\p \cc$ and $\p \cc'$
exist.

\begin{proof}
First suppose that $R$ takes $F'$ into $F$ and takes 
$\inj C'$ into $\inj C$.
Since $R$ commutes
with finite limits, it takes
special $F'$-maps to special $F$-maps.  This means that $R$ takes retracts of
special $F'$-maps to retracts of special $F$-maps, which means that
$R$ preserves fibrations.

To show that $R$ preserves acyclic fibrations, recall that an
acyclic fibration in $\p \cc'$ is a retract of a special $\inj
C'$-map (see Proposition \ref{pro:a-fib}).  We can use the same
argument as in the previous paragraph, using that $R$ takes $\inj
C'$ to $\inj C$.

Now suppose that $L$ and $R$ are a Quillen adjoint pair.
Let $p:X \map Y$ belong to $F'$.  We want to show that $Rp$ belongs
to $F$.  The constant pro-map $cp: cX \map cY$ is a special $F'$-map,
so it is a fibration in $\p \cc'$.  Since $R$ is a right Quillen functor,
$R(cp)$ must be a fibration in $\p \cc$, so it is a retract of a special
$F$-map $q: Z \map W$.
Note that $R(cp)$ equals $c(Rp)$.

We have a diagram
\[
\xymatrix{
cRX \ar[r] \ar[d] & Z \ar[r]\ar[d] & cRX \ar[d] \\
cRY \ar[r]  & W \ar[r] & cRY                    }
\]
in $\p \cc$, which is represented by a diagram
\[
\xymatrix{
RX \ar[r] \ar[d] & Z_s \ar[r]\ar[d] & RX \ar[d] \\
RY \ar[r]  & W_s \ar[r] & RY                    }
\]
in $\cc$.  By  Lemma \ref{lem:fib} each $Z_s \map W_s$ belongs to
$F$. Now $F$ is closed under retract by Axiom \ref{ax:retract}, so
$RX \map RY$ belongs to $F$.  This is what we wanted to prove.

Using Proposition \ref{pro:a-fib},
the proof that $R(\inj C')$ is contained in $R(\inj C)$ is identical.
Note that $\inj C$ is closed under finite compositions and arbitrary
base changes for formal reasons.
\end{proof}

Although Theorem \ref{thm:Q-adj} gives elegant necessary and sufficient
conditions for a Quillen adjunction, it can sometimes be hard to verify
these conditions in practice.  We give other conditions that 
are sometimes easier to check.

\begin{pro}
\label{prop:Q-adj1} Let $(A, C, W, F)$ and $(A', C', W', F')$ be
proper filtered model structures on the categories $\cc$ and $\cc'$
respectively. Let $L : \cc \rarr \cc'$ be a left adjoint of $R :
\cc' \rarr \cc$. Suppose that for every $b$ in $A'$, there exists
$a$ in $A$ such that $L(C_a)$ is contained in $C'_{b}$ and $L(W_a
\cap C_a)$ is contained in $W'_{b} \cap C'_{b}$. Then the induced
functors $L : \p \cc \rarr \p \cc'$ and $R : \p \cc' \rarr \p \cc$
are a Quillen adjoint pair when $\p \cc$ and $\p \cc'$ are equipped
with the model structures of Theorem \ref{thm:ms}.
\end{pro}

\begin{proof}
To see that $L$ preserves cofibrations, use 
Corollary \ref{cor:preserve} and the definition of cofibrations.

Now consider an acyclic cofibration $i$.
By Proposition
\ref{pro:a-cofib}, we know that $i$ is an essentially levelwise
$\proj F_a$-map for every $a$.  
Using Axiom \ref{ax:acyclic}, it
follows that $i$ is an essentially levelwise $(W_a \cap C_a)$-map
for every $a$. Now apply Corollary \ref{cor:preserve}.
\end{proof}

We now consider Quillen equivalences.

\begin{thm}
\label{thm:Q-eq}
Let $(A, C, W, F)$ and $(A', C', W', F')$ be
proper filtered model structures on the categories $\cc$ and $\cc'$
respectively. Let $L : \cc \rarr \cc'$ be a left adjoint of $R :
\cc' \rarr \cc$ such that
the induced functors $L: \p \cc \rarr \p \cc'$ and 
$R : \p \cc' \rarr \p \cc$ are a Quillen
adjoint pair when $\p \cc$ and $\p \cc'$ are equipped with the model
structures of Theorem \ref{thm:ms}.
Suppose also that:
\begin{enumerate}
\item
For every $b$ in $A'$, there exists $a$ in $A$ such that if
$X \map RY$ is a map in $\cc$ belonging to $W_a$, $* \map X$ belongs
to some $C_c$, and $Y \map *$ belongs to $F'$; then the adjoint
map $LX \map Y$ belongs to $W'_b$.
\item
For every $b$ in $A$, there exists $a$ in $A'$ such that if
$LX \map Y$ is a map in $\cc'$ belonging to $W'_a$, $* \map X$ belongs
to some $C_c$, and $Y \map *$ belongs to $F'$; then the adjoint
map $X \map RY$ belongs to $W_b$.
\end{enumerate}
Then $L$ and $R$ are a Quillen equivalence between $\p \cc$
and $\p \cc'$.
\end{thm}

\begin{proof}
Suppose that $X$ is a cofibrant object of $\p \cc$ and $Y$ is a fibrant
object of $\p \cc'$.
Since $X$ is cofibrant,
we may assume
that each map $* \map X_s$ belongs to some $C_c$.
Lemma \ref{lem:level-fibrant} says that we may assume that each
map $Y_t \map *$ belongs to $F'$.

Now suppose that $f:LX \map Y$ is a weak equivalence in $\p \cc'$.
Our goal is to show that the adjoint map $g:X \map RY$
is a weak equivalence in $\p \cc$.  
Using the level replacement of
\cite[App.~3.2]{AM}, we may reindex $X$ and $Y$ in such a way that
$f$ is a level map, each $* \map X_s$ still belongs to some $C_c$,
and each $Y_t \map *$ still belongs to $F'$.

Factor $f$ into an acyclic cofibration $LX \map Z$
followed by a fibration $Z \map Y$.  If we use the method
of Lemma \ref{lem:factor2} to produce this factorization, we find 
that for all $a$ in $A'$, there exists an $s(a)$ such that
the map $LX_t \map Z_t$ belongs to $W_a$ for all $t \geq s(a)$.
Moreover, $Z \map Y$ is a levelwise $F'$-map by Lemma \ref{lem:fib}, so 
$Z_t \map *$ belongs to $F'$ for all $t$.

Now condition (2) implies that for every $b$ in $A$,
there exists $a$ in $A'$ such that 
the map $X_t \map RZ_t$ belongs to $W_b$ for all $t \geq s(a)$.
In particular, the map
$X \map RZ$ is an essentially levelwise $W_b$-map for every $b$.
Thus $X \map RZ$ is a weak equivalence.

The map $g$ factors as
$X \map RZ \map RY$.
We have just observed that the first map is a weak equivalence.
For the second map, note that the two-out-of-three axiom implies that
$Z \map Y$ is an acyclic fibration and that the
right Quillen functor $R$ preserves acyclic fibrations.
This shows that $g$ is a weak equivalence and finishes one half
of the proof.

Now suppose that $g:X \map RY$
is a weak equivalence in $\p \cc$.
Our goal is to show that the adjoint map 
$f:LX \map Y$ 
is a weak equivalence in $\p \cc'$.
Using the level replacement of
\cite[App.~3.2]{AM}, we may reindex $X$ and $Y$ in such a way that
$f$ is a level map, each $* \map X_s$ still belongs to some $C_c$,
and each $Y_t \map *$ still belongs to $F'$.

Use the method
of Lemma \ref{lem:factor2} to 
factor $g$ into a cofibration $X \map Z$
followed by an acyclic fibration $Z \map RY$.  
In the notation of the proof of that lemma, we may choose
$\phi(t)$ sufficiently large such that $X_t \map Z_t$ and
$* \map X_t$ both belong to $C_c$ for some value of $c$.
Then the composition $* \map Z_t$ also belongs to $C_c$.

Also, $Z \map RY$ is a special $\inj C$-map.
The proof of Lemma \ref{lem:fib} can be adapted line by line to show
that $Z \map RY$ is a levelwise $\inj C$-map.  Lemma \ref{lem:inj-C}
implies that $Z \map RY$ is a levelwise $W_a$-map for every $a$.
Now condition (1) of the theorem implies that the adjoint map
$LZ \map Y$ is a levelwise $W'_b$-equivalence for all $b$ in $A'$.

The map $f$ factors as
\[
\xymatrix@1{
LX \ar[r] & LZ \ar[r] & Y.   }
\]
We have just observed that the second map is a weak equivalence.
The two-out-of-three axiom implies that $X \map Z$ is an
acyclic cofibration.  Since the left Quillen functor $L$
preserves acyclic cofibrations,
it follows that $LX \map LZ$ is also a
weak equivalence.  Thus $f$ is a weak equivalence.
\end{proof}

We can apply our general results above to the specific case of 
strict model structures.  It was an oversight that this result did
not appear in \cite{strict}.

\begin{thm}
\label{thm:strict-Q}
Let $L: \cc \map \cc'$ and $R: \cc' \map \cc$ be a Quillen pair between
two proper model categories $\cc$ and $\cc'$.
The induced functors $L: \p \cc \rarr \p \cc'$ and 
$R : \p \cc' \rarr \p \cc$ are a Quillen
adjoint pair when $\p \cc$ and $\p \cc'$ are equipped with
strict model structures.
If $L$ and $R$ are a Quillen equivalence, then 
$L$ and $R$ induce a Quillen equivalence between 
the strict model structures on $\p \cc$ and $\p \cc'$.
\end{thm}

\begin{proof}
Recall from Proposition \ref{exmp:strict} that a strict model structure
is the model structure associated to a filtered model structure indexed
by a set with only one element.
In the case when $A$ and $A'$ have only one element, the conditions 
of Theorem \ref{thm:Q-adj} reduce to the definition of a Quillen pair.
Similarly, the conditions of Theorem \ref{thm:Q-eq} reduce to the 
definition of a Quillen equivalence.
\end{proof}

\section{Examples}
\label{exmp-gen}

In this section, we describe several examples of model structures on
pro-categories that can be established with Theorem \ref{thm:ms}.

\begin{exmp}
We give a concrete example that is an application of Proposition
\ref{pro:actual-mc}. Let $\cc$ be the category of $S$-modules
\cite{ekmm}, and let $A = \mathbb{N}$. For each $n$, let $E_n$ be a
generalized homology theory such that if a map of spectra is an
$E_n$-homology equivalence, then it is also an $E_m$-homology
equivalence for all $m \leq n$. Let $W_n$ be the class of
$E_n$-homology equivalences. For each $n$ there is a proper model
structure on $\cc$ whose weak equivalences are the class $W_n$ and
whose cofibrations are the class $C$ of all cofibrations of spectra
\cite[VIII.1]{ekmm}.
Let $C_n = C$ for all $n$. Now the classes $C_n$ satisfy the
necessary inclusion relationships trivially, and the classes $W_n$
satisfy the necessary inclusion relationships by assumption on the
homology theories $E_n$. By Proposition \ref{pro:actual-mc}, there
is a model structure on $\p \cc$ such that the weak equivalences are
maps that are \ess\ $E_n$-homology equivalences for all $n$. If
$L_n$ is the $E_n$-localization functor, then the map $cX \rarr \{
L_n X \}$ is a weak equivalence for all spectra $X$. Note that there
are natural transformations $L_n \rarr L_{n-1}$ because of the
assumption on the homology theories $E_n$.

The role of $S$-modules in this example is not central.  Any of the
standard model structures for stable homotopy theory will work 
just as well.
\end{exmp}

\begin{exmp}
\label{ex:space} We give a concrete illustration of Proposition
\ref{pr:cof-gen}. Let $\cc$ be the category of simplicial sets. Let
$A$ be $\nn$, and define $I_n$ to be the set of generating
cofibrations $\partial \Delta^k \map \Delta^k$ with $k > n$. The
results of \cite[Sec.~3]{pro} imply that $W_n$ is the class of
$n$-equivalences, and $F_n$ is the class of co-$n$-fibrations, i.e.,
the class of maps that are both fibrations and co-$n$-equivalences.
The hypotheses of Proposition \ref{pr:cof-gen} are satisfied
because the class of $n$-equivalences is closed under retract and
because if any two of $f$, $g$, and $gf$ are $n$-equivalences, then
the third is an $(n-1)$-equivalence.

Before we can use Theorem \ref{thm:ms} to obtain a model structure on
pro-simplicial sets, we must also observe that the filtered model structure of
the previous paragraph is proper.  This follows from the results of
\cite[Sec.~3]{pro}; see the end of the proof of
Proposition \ref{prop:equiv-proper} and Lemma \ref{lem:cobase} for the
topological analogue.

The resulting model structure on the category of pro-simplicial sets
is the same as the model structure of \cite[Thm.~6.4]{pro}. The weak
equivalences are the pro-maps $f$ such that for every $n \geq 0$,
$f$ is isomorphic to a levelwise $n$-equivalence. This fact is not
clearly stated in \cite{pro}, but see \cite[Prop.~6.8]{pro} for a
``morally equivalent'' claim. Compared to the technical methods of
\cite{pro} involving basepoints, this approach is much simpler.
\end{exmp}

\begin{exmp}
\label{ex:top} Similarly to Example \ref{ex:space}, there is a
filtered model structure on the category of topological spaces.  The
class $C_n$ is the class of Serre cofibrations, i.e., retracts of
relative cell complexes; the class $W_n$ is the class of
$n$-equivalences; and the class $F_n$ is the class of
co-$n$-fibrations, i.e., Serre fibrations that are also
co-$n$-equivalences.  In this context, $I_n$ is the set of
generating cofibrations $S^{k-1} \map D^k$ with $k > n$.
In the same way as the previous example, one can show that
the hypotheses of Proposition \ref{pr:cof-gen} are satisfied and that
the associated filtered model structure is proper.  The resulting
model structure
on pro-spaces is Quillen equivalent to the model structure on pro-simplicial
sets from the previous example.
\end{exmp}

\begin{exmp}
\label{ex:spectra} The following example is a stable version of
Example \ref{ex:space}. Let $\cc$ be the category of spectra.  Let
$A$ be $\zz$, and define $I_n$ to be the set of generating
cofibrations whose cofibers are spheres of dimension greater than
$n$.  The results of \cite[Sec.~4]{ipi} imply that $W_n$ is the
class of $n$-equivalences, and $F_n$ is again the class of
co-$n$-fibrations. As before, the hypotheses of Proposition \ref{pr:cof-gen}
are satisfied, and the associated filtered
model structure is proper. The resulting model structure
on the category of pro-spectra is the same as the model structure of
\cite{ipi}. The weak equivalences can be described in terms of
pro-homotopy groups, but the reformulation is not quite as obvious
as one might expect.  See \cite{ipi} for details.
\end{exmp}

\begin{exmp} \label{hpostnikov}
Let $\cc$ be the category of spaces, and let $A$ equal $\nn$. Let
$h_*$ be a homology theory on $\cc$ that satisfies the colimit
axiom.  Let $W_n$ be the class of maps $f$ such that $h_i ( f)$ is
an isomorphism for $i < n $ and $ h_n (f)$ is a surjection. In order
to obtain a {\em proper} filtered model structure, we must assume
that $W_n$ is preserved by base changes along fibrations. If we let
$C_n$ be the class of cofibrations and $F_n$ be the class $\inj (
W_n \cap C_n )$, then we get a proper filtered model structure on
spaces.

We outline the verification of the axioms for a filtered model
structure in this case.  Axioms \ref{ax:2/3} and \ref{ax:retract}
are obvious, as is the first half of Axiom \ref{ax:acyclic}.  We
defer the second half of Axiom \ref{ax:acyclic} until later.

The first factorization of Axiom \ref{ax:factor} is given by
factorizations into cofibrations followed by acyclic fibrations. For
the second part of this axiom, an adaptation of the small object
argument in \cite[Sec.~11]{bou} gives the desired factorization when
the source is cofibrant. The basic idea is to replace all statements
of the form ``$h_*(K,L) = 0$'' with statements of the form
``$h_i(K,L) = 0$ for $i \leq n$''. The factorization in
general follows from standard arguments with model categories, together
with properness for the category of spaces.

Now we return to the second half of Axiom \ref{ax:acyclic}. Using
the factorization of the previous paragraph, a retract argument
shows that $\proj F_n$ equals $W_n \cap C_n$.

Having identified $\proj F_n$, the factorizations required by Axiom
\ref{ax:lift} are provided by factorizations into cofibrations
followed by acyclic fibrations.

Axiom \ref{ax:cobase} follows from consideration of the long exact
sequence in homology associated to a cofiber sequence. Axiom
\ref{ax:base} is satisfied by our assumption on $W_n$.

If $h_*$ is a periodic cohomology theory, then the model structure
of this example is just the strict model structure associated to the
$h_*$-local model structure on spaces \cite{bou}.
\end{exmp}

\section{The underlying model structure for pro-finite groups}
\label{section-free}

Recall our goal of finding a model structure for pro-$G$-spaces in
which $\{ E(G/U) \}$ is a cofibrant replacement for $*$, where $G$
is a pro-finite group and $U$ ranges over all finite-index normal
subgroups of $G$. We begin by describing the equivariant analogue of
Example \ref{ex:space}, but this model structure turns out not to
have the desired property.

Let $G$ be a (topological) group, and let $\cc$ be the category of
$G$-spaces. Let $C$ be the class of retracts of relative $G$-cell
complexes; this is the class of cofibrations in one of the usual
model structures on the category of $G$-spaces. Now let $A$ be $\nn$,
and let $I_n$ be the set of generating cofibrations of the form
$S^{k-1} \times G/H \map D^k \times G/H$, where $k > n$ and $H$ is a
(closed) subgroup of $G$. The framework of Proposition \ref{pr:cof-gen}
applies; the hypotheses of this result can be
verified just as in \cite[Sec.~3]{pro}. The class $W_n$ turns out to
be the class of $G$-equivariant $n$-equivalences, i.e., the maps $f:
X \map Y$ such that $f^H: X^H \map Y^H$ is an $n$-equivalence for
all (closed) subgroups $H$ of $G$. The class $F_n$ is the class of
equivariant co-$n$-fibrations, i.e., the maps $f: X \map Y$ such
that $f^H: X^H \map Y^H$ is a co-$n$-fibration for all (closed)
subgroups $H$ of $G$.

The resulting model structure is a $G$-equivariant analogue of the
model structure of \cite[Thm.~6.4]{pro}. It can be shown that a map
$f: X \map Y$ of pro-$G$-spaces is a weak equivalence if and only if
$\underline{\pi}_n f: \underline{\pi}_n X \map \underline{\pi}_n Y$
is a pro-isomorphism of pro-coefficient systems  for all $n \geq 0$.

The map $\{ E(G/U) \} \rarr c (* )$ induces a pro-isomorphism after
applying $\pi_n^V$ for any fixed $V$.  However, it does {\em not}
induce a pro-isomorphism after applying $\underline{\pi}_n$.  The
problem is that in choosing refinements, one must make different
choices for different values of $V$.  There is no one choice that
works for all $V$.

The point of the previous paragraphs is that we must work harder to
obtain our desired model structure.  The rest of this section
provides the details.

Let $G$ be a pro-finite group.  This means that $G$ is a topological
group such that $G \rarr \lim_U G/U$ is an isomorphism, where $U$
ranges over all normal subgroups of $G$ such that $G/U$ is discrete
and finite. Usually, a pro-finite group is  viewed as a topological
group, where the topology is totally disconnected, compact, and
Hausdorff. It is also possible to think of $G$ as a pro-object in
the category of finite groups. We will use both viewpoints.

Let $\cc$ be the category of compactly generated weak Hausdorff
spaces equipped with a continuous $G$-action. First, recall that
$\cc$ is a simplicial category.  If $X$ is a $G$-space and $K$ is a
simplicial set, then $X \otimes K$ is defined to be $X \times |K|$,
where the realization $|K|$ has a trivial $G$-action. Also $X^K$ is
the topological mapping space $F(|K|, X)$ of non-equivariant
continuous maps $|K| \map X$.  If $X$ and $Y$ are both $G$-spaces,
then $\Map(X,Y)$ is the simplicial set whose $n$-simplices are
equivariant maps $X \times |\Delta^n| \map Y$.

\begin{defn}
\label{defn:U} Let $U$ be a finite-index subgroup of $G$.  A
$G$-equivariant map $f: X \map Y$ is a \mdfn{$U$-weak equivalence}
if the map $f^V$ on $V$-fixed points is a weak equivalence for all
finite-index subgroups $V$ of $U$, and it is a \mdfn{$U$-fibration}
if $f^V$ is a fibration for all finite-index subgroups $V$ of $U$.
It is a \mdfn{$U$-cofibration} if it is a retract of a relative cell
complex built from cells of the form $G/V \times S^n \map G/V \times
D^{n+1}$, where $V$ is a finite-index subgroup of $U$.
\end{defn}

Recall that an $n$-cofibration \cite[Defn.~3.2]{pro} is a
cofibration that is also an $n$-equivalence (i.e., an isomorphism on
homotopy groups up to dimension $n-1$ and a surjection on $n$th
homotopy groups). Dually, a co-$n$-fibration is a fibration that is
also a co-$n$-equivalence (i.e., an isomorphism on homotopy groups
above dimension $n$ and an injection on $n$th homotopy groups). We
will now make similar definitions in the equivariant situation.

\begin{defn}
\label{defn:U-n-equiv}
A map $f$ is a \mdfn{$U$-$n$-equivalence} (resp.,
\mdfn{$U$-co-$n$-equivalence} if the map $f^V$ is an $n$-equivalence
(resp., co-$n$-equivalence) for all finite-index subgroups $V$ of
$U$.  A map $f$ is a \mdfn{$U$-co-$n$-fibration} if $f^V$ is a
co-$n$-fibration for all finite-index subgroups $V$ of $U$. Finally,
a map is a \mdfn{$U$-$n$-cofibration} if it is a $U$-cofibration
such that $f^V$ is an $n$-equivalence for all finite-index subgroups
$V$ of $U$.
\end{defn}

The following two lemmas are proved in exactly the same way as their
non-equivariant analogues \cite[Sec.~3]{pro}.

\begin{lem}
\label{lem:equiv-factor} Any $G$-equivariant map can be factored
into a $U$-$n$-cofibration followed by a $U$-co-$n$-fibration.
\end{lem}

\begin{proof}
Use the small object argument applied to the set of maps of the form
$G/V \times S^k \map G/V \times D^{k+1}$, where $V$ is a
finite-index subgroup of $U$ and $k \geq n$, together with all maps
of the form $G/V \times I^m \map G/V \times I^{m+1}$, where $V$ is a
finite-index subgroup of $U$ and $m$ is arbitrary.  Here $I^m$ is
the $m$-cube, and the map $I^m \map I^{m+1}$ is the inclusion of a
face.
\end{proof}

\begin{lem}
\label{lem:equiv-lift} The classes of $U$-$n$-cofibrations and
$U$-co-$n$-fibrations are determined by lifting properties with
respect to each other.
\end{lem}

\begin{proof}
This can be proved with an obstruction theory argument.  See
\cite[Lem.~3.4]{pro} and \cite[Lem.~3.6]{pro} for more details.
\end{proof}

Let \mdfn{$A$} be the set consisting of all pairs $(U,n)$, where $U$ is a
finite-index subgroup of $G$ and $n$ is a non-negative integer.  We
write $(U,n) \geq (V,m)$ if $U$ is contained in $V$ and if $n \geq
m$.  This makes $A$ into a directed set.  In other words, given
$(U_1, n_1)$ and $(U_2, n_2)$, there exists $(V,m)$ such that $(V,m)
\geq (U_1, n_1)$ and $(V,m) \geq (U_2, n_2)$.  To see why this is
true, just observe that $U_1 \cap U_2$ is a finite-index subgroup of
$G$ whenever $U_1$ and $U_2$ are finite-index subgroups.

For each $(U,n)$ in $A$, we will define three classes $C_{U,n}$,
$W_{U,n}$, and $F_{U,n}$ of $G$-equivariant maps.

\begin{defn}
\label{defn:CWF-equivariant} The class \mdfn{$C_{U,n}$} is the class
of all $U$-cofibrations.  The class \mdfn{$W_{U,n}$} is the class of
maps that are $V$-$n$-equivalences for some $V$.  The class
\mdfn{$F_{U,n}$} is the class of all $U$-co-$n$-fibrations.
\end{defn}

Note that $C_{U,n}$ does not actually depend on $n$, and $W_{U,n}$
does not actually depend on $U$.  The point of this seemingly
confusing notation is that there is just one indexing set $A$ for
all three families of classes.

As in Definition \ref{defn:F-inj-C}, we write \mdfn{$F$} for the
union of the classes $F_{U,n}$.

\begin{exmp}
The following example emphasizes a subtlety in the definition of $W_{U,n}$,
Consider the map
\[
\amalg_U E(G/U) \rarr \amalg_U *,
\]
where $U$ ranges over the finite-index normal subgroups of $G$. This
map is an underlying weak equivalence.  However, it is not a
$V$-equivalence for any finite-index subgroup $V$ of $G$ and thus does
not belong to $W_{V,n}$ for any $(V,n)$.
\end{exmp}

\begin{lem}
\label{lem:contain-equivariant} If $(V,m) \geq (U,n)$, then
$C_{V,m}$ is contained in $C_{U,n}$, $W_{V,m}$ is contained in
$W_{U,n}$, and $F_{U,n}$ is contained in $F_{V,m}$.
\end{lem}

\begin{proof}
All three claims follow immediately from the definitions.  If $V$ is
contained in $U$, then the set of generating $V$-cofibrations is a subset of
the set generating $U$-cofibrations. This shows that $C_{V,m}$ is
contained in $C_{U,n}$.

If $m \geq n$, then an $m$-equivalence is automatically an
$n$-equivalence; this shows that $W_{V,m}$ is contained in
$W_{U,n}$.

If $m \geq n$, then a co-$n$-fibration is automatically a
co-$m$-fibration. If $V$ is contained in $U$, then a $U$-fibration is
automatically a $V$-fibration. This shows that $F_{U,n}$ is
contained in $F_{V,m}$.
\end{proof}

\begin{lem}
\label{lem:equivariant-inj-C-F} The class $\inj C_{U,n}$ (i.e., the
class of maps that have the right lifting property with respect to
all $U$-cofibrations equals the class of $U$-acyclic fibrations
(i.e., maps that are both $U$-weak equivalences and $U$-fibrations).
The class $\proj F_{U,n}$ (i.e., the class of maps that have the
left lifting property with respect to all $U$-co-$n$-fibrations) equals
the class of $U$-$n$-cofibrations.
\end{lem}

\begin{proof}
The first claim follows from standard equivariant homotopy theory.
The second claim is immediate from Lemma \ref{lem:equiv-lift}.
\end{proof}

Recall that $\cc$ is the category of compactly generated weak
Hausdorff spaces with continuous $G$-actions and $G$-equivariant
maps. The following definition is a special case of Definitions
\ref{defn:cofib}, \ref{defn:we}, and \ref{defn:fib}.

\begin{defn}
\label{defn:equiv-ms} A map in $\p \cc$ is a \mdfn{cofibration} if
it is an essentially levelwise $C_{U,n}$-map for every $(U,n)$ in
$A$. A map in $\p \cc$ is a \mdfn{weak equivalence} if it is an
essentially levelwise $W_{U,n}$-map for every $(U,n)$ in $A$. A map
in $\p \cc$ is a \mdfn{fibration} if it is a retract of a special
$F$-map.
\end{defn}

\begin{thm}
\label{thm:equiv-ms} Definition \ref{defn:equiv-ms} is a proper
simplicial model structure on the category $\p \cc$.
\end{thm}

\begin{proof}
Using Theorems \ref{thm:ms} and \ref{thm:simplicial}, we just need
to verify that Definition \ref{defn:CWF-equivariant} is a proper
simplicial filtered model structure.  This is provided below in
Propositions \ref{prop:equiv-filtered}, \ref{prop:equiv-proper}, and
\ref{prop:equiv-simplicial}.
\end{proof}

\begin{pro}
\label{prop:equiv-filtered} Definition \ref{defn:CWF-equivariant} is
a filtered model structure.
\end{pro}

\begin{proof}
We have to verify Axioms \ref{ax:2/3} through \ref{ax:lift}.

We have already observed that $A$ is a directed set. Lemma
\ref{lem:contain-equivariant} says that the containments given in
Definition \ref{defn:filtered} are satisfied. Also, the category of
$G$-spaces is complete and cocomplete; limits and colimits are
constructed in the underlying category of topological spaces.

For Axiom \ref{ax:2/3}, first observe that if any two of the maps
$f$, $g$, and $gf$ are $V$-$(n+1)$-equivalences, then a simple
diagram chase shows that the third is a $V$-$n$-equivalence.  If any
two of $f$, $g$, and $gf$ belong to $W_{U,n+1}$, then there exists a
finite-index subgroup $V$ of $G$ such that the two maps are
$V$-$(n+1)$-equivalences.  The third map is a $V$-$n$-equivalence,
which means that it belongs to $W_{U,n}$.

We now consider Axiom \ref{ax:retract}.  The 
$V$-$n$-equivalences are closed under retract since  non-equivariant
$n$-equivalences are closed under retract, so $W_{U,n}$ is closed
under retract.
The class $C_{U,n}$ of $U$-cofibrations is closed under retract, finite
compositions, and  arbitrary cobase changes because it is defined
in terms of retracts of relative cell complexes.
The  class $F_{U,n}$ is defined by a right lifting property (see
Lemma \ref{lem:equiv-lift}), so it is closed under retract,
finite compositions, and arbitrary base changes.

Axiom \ref{ax:acyclic} is immediate from Lemma
\ref{lem:equivariant-inj-C-F}.

The first half of Axiom \ref{ax:factor} is given by factorizations
into $U$-cofibrations followed by $U$-acyclic fibrations; these
factorizations are supplied by standard equivariant homotopy theory.
The second half of Axiom \ref{ax:factor} is given by Lemma
\ref{lem:equiv-factor}.

For Axiom \ref{ax:lift}, let $f$ belong to $W_{U,n}$.  This means
that $f$ is a $V$-$n$-equivalence for some $V$.  Now factor $f$ into
a $V$-cofibration $i$ followed by a $V$-acyclic fibration $p$. By
Lemma \ref{lem:equivariant-inj-C-F}, this means that $p$ belongs to
$\inj C_{V,n}$.
Because $f$ is a $V$-$n$-equivalence, $i$ is also a
$V$-$n$-equivalence and hence a $V$-$n$-cofibration.  By Lemma
\ref{lem:equivariant-inj-C-F} again, this means that $i$ belongs to
$\proj F_{V,n}$.
\end{proof}

\begin{pro}
\label{prop:equiv-proper} Definition \ref{defn:CWF-equivariant} is a
proper filtered model structure.
\end{pro}

\begin{proof}
We just need to verify Axioms \ref{ax:cobase} and \ref{ax:base}.

For Axiom \ref{ax:cobase}, suppose that $f: A \map B$ is a
$U$-cofibration and that $g: A \map C$ is a $V$-$n$-equivalence for
some finite-index subgroups $U$ and $V$ of $G$. We may replace $V$
by $V \cap U$ to assume that $V$ is a finite-index subgroup of $U$;
this is allowed because $g$ is a still a $V$-$n$-equivalence. We
will show that the map $h: B \map B \amalg_A C$ is also a
$V$-$n$-equivalence. If $W$ is a finite-index subgroup of $V$, then
$(B \amalg_A C)^W$ is equal to $B^W \amalg_{A^W} C^W$; this uses the
fact that $f$ is injective. Now $W$ is a finite-index subgroup of
$U$, so $f^W$ is a non-equivariant cofibration. Since $g^W$ is a
non-equivariant $n$-equivalence, we need only show that cobase
changes along non-equivariant cofibrations preserve non-equivariant
$n$-equivalences.  This is proved in Lemma \ref{lem:cobase} below.

For Axiom \ref{ax:base}, suppose that $f:X \map Y$ is a
$U$-fibration and that $g: Z \map Y$ is a $V$-$n$-equivalence for
some finite-index subgroups $U$ and $V$ of $G$. Choose a
finite-index subgroup $W$ contained in both $V$ and $U$.  Then $f$
is a $W$-fibration and $g$ is a $W$-$n$-equivalence.  Taking fixed
points commutes with fiber products.  Therefore, in order to show
that $X \times_Y Z \map X$ is a $W$-$n$-equivalence, we only need
prove that base changes of non-equivariant $n$-equivalences along
non-equivariant fibrations are $n$-equivalences.  This last fact
follows from the five lemma and the long exact sequence of homotopy
groups for a fibration.
\end{proof}

\begin{lem}
\label{lem:cobase} Let $f: A \map B$ be a cofibration of topological
spaces, and let $g: A \map C$ be an $n$-equivalence.  Then the map
$h:B \map B \amalg_A C$ is also an $n$-equivalence.
\end{lem}

\begin{proof}
Since the usual model category on topological spaces is left proper,
the pushout $B \amalg_A C$ is in fact a homotopy pushout. Therefore,
we may replace $g$ by a weakly equivalent cofibration; this will not
change the weak homotopy type of $h$.

Now we have that $g$ is an $n$-cofibration (i.e., a cofibration and
an $n$-equivalence). The class of $n$-cofibrations is determined by
a left lifting property; this is the non-equivariant version
of Lemma \ref{lem:equiv-lift}.
Therefore, the class of $n$-cofibrations is
closed under arbitrary cobase changes, so $h$ is also an
$n$-cofibration.
\end{proof}

\begin{rem}
The reader may feel that it is not possible to prove Lemma
\ref{lem:cobase} without using van Kampen's theorem.  Van Kampen's
theorem is necessary to prove that the model category of topological
spaces is left proper, so we are in fact using it in a disguised
way.
\end{rem}

\begin{pro}
\label{prop:equiv-simplicial} Definition \ref{defn:CWF-equivariant}
is a simplicial filtered model structure.
\end{pro}

\begin{proof}
We have already observed that $\cc$ is a simplicial category, so we
just need to prove that Axiom \ref{ax:simplicial} holds.

Let $j: K \map L$ be a cofibration of finite simplicial sets, and
let $i: A \map B$ be a $U$-cofibration.  Standard equivariant
homotopy theory implies that the map
\[
f: A \otimes L \amalg_{A \otimes K} B \otimes K \map B \otimes L
\]
is also a $U$-cofibration.

Similarly, if $j$ is an acyclic cofibration, then standard
equivariant homotopy theory implies that $f$ is a $U$-acyclic
cofibration.  This implies that $f$ is a $U$-$n$-cofibration, which
means that it belongs to $\proj F_{U,n}$ by
Lemma \ref{lem:equivariant-inj-C-F}.

Next, suppose that $j$ is a cofibration and that $i$ belongs to
$\proj F_{U,n}$.  By Lemma \ref{lem:equivariant-inj-C-F}, this means
that $i$ is a $U$-$n$-cofibration. We want to conclude that $f$ is
also a $U$-$n$-cofibration.  We have already shown that $f$ is a
$U$-cofibration, so we just need to show that $f^V$ is an
$n$-equivalence for every finite-index subgroup of $V$ of $U$.

The map $f^V$ is equal to the map
\[
A^V \otimes L \amalg_{A^V \otimes K} B^V \otimes K \map B^V \otimes
L.
\]
This follows from the fact that the $G$-actions on $K$ and $L$ are
trivial and that taking $V$-fixed points commutes with the pushout
because $A \otimes K \map A \otimes L$ is injective. Now the map
$i^V: A^V \map B^V$ is an $n$-equivalence because $i$ is a
$U$-$n$-cofibration, so the desired conclusion follows from Lemma
\ref{lem:simplicial} below.
\end{proof}

\begin{lem}
\label{lem:simplicial} Suppose that $j: K \map L$ is a cofibration
of simplicial sets, and suppose that $i: A \map B$ is an
$n$-cofibration of non-equivariant topological spaces.  Then the map
\[
f: A \otimes L \amalg_{A \otimes K} B \otimes K \map B \otimes L
\]
is an $n$-cofibration.
\end{lem}

\begin{proof}
Consider the diagram
\[
\xymatrix{
A \otimes K \ar[r]\ar[d] & B \otimes K \ar[d] \\
A \otimes L \ar[r]\ar[drr] &
  A \otimes L \amalg_{A \otimes K} B \otimes K \ar[dr]^f \\
& & B \otimes L.  }
\]
We will show that $A \otimes K \map B \otimes K$ and $A \otimes L
\map B \otimes L$ are $n$-cofibrations.  Then, since
$n$-cofibrations are preserved by cobase changes, it follows that
the map $A \otimes L \map A \otimes L \amalg_{A \otimes K} B \otimes
K$ is an $n$-cofibration.  A small diagram chase proves that $f$ is
an $n$-equivalence.

It remains only to prove that $A \otimes K \map B \otimes K$ is an
$n$-cofibration; the proof for $A \otimes L \map B \otimes L$ is
identical. First, standard homotopy theory of topological spaces
says that $A \times |K| \map B \times |K|$ is a cofibration. Second,
since homotopy groups commute with products, the map $A \times |K|
\map B \times |K|$ is an $n$-equivalence.
\end{proof}

\begin{rem}
\label{rem:stable-equiv} The basic ideas of this section can be
implemented in exactly the same way for naive $G$-spectra.  One
small difference is that the indexing set $A$ consists of pairs
$(U,n)$ where $U$ is a finite-index subgroup of $G$ as before but
$n$ is an arbitrary integer, possibly negative. See \cite{ipi} for
details concerning $n$-cofibrations and co-$n$-fibrations of
spectra.
\end{rem}

Finally, we can establish our main motivation for producing the
model structure of Theorem \ref{thm:equiv-ms}.

\begin{pro}
\label{pro:equiv-cofib} The pro-$G$-space $\{ E(G/U) \}$ is a
cofibrant replacement for the constant trivial pro-space $c(*)$ in
the model structure of Theorem \ref{thm:equiv-ms}.
\end{pro}

\begin{proof}
Note first that for each finite-index subgroup $U$, $E(G/U)$
can be built from cells of the form $S^{k-1} \times G/U \map D^k
\times G/U$. This means that the map from the empty set to $E(G/V)$
is a $U$-cofibration whenever $V$ is contained in $U$.
It follows that $\phi \map \{ E(G/U) \}$ is a
cofibration of pro-$G$-spaces.

Now we will show that the map $E(G/U) \map *$ is a $U$-weak
equivalence for each $U$; this will imply that the map $\{E(G/U) \}
\map c(*)$ is a weak equivalence of pro-$G$-spaces. If $V$ is a
finite-index subgroup of $U$ then the $V$-fixed points of $E(G/U)$
equals $E(G/U)$. Since $E(G/U)$ is contractible, it follows that the
map
\[
E(G/U)^V \map *^V
\]
is a weak equivalence.
\end{proof}

\end{document}